\documentclass{mcom-l}

\usepackage{amssymb,amsthm,amsfonts,amsmath,graphicx,array,url,subfigure}

\newcounter{lineno}
\newcounter{linelineno}[lineno]
\newcounter{linelinelineno}[linelineno]

\newcounter{algnum}
\renewcommand{\thealgnum}{\arabic{algnum}}

\newenvironment{algorithm}[3]
{
        \refstepcounter{algnum}
        \bigskip\goodbreak\hrule\medskip
                \leftline{{\sc Algorithm \thealgnum:} \bf #1}
                \medskip\hrule\medskip\nobreak
                \tt
        \begin{tabbing}
        OUTPUT: \= \kill
        INPUT:  \>#2 \\
        OUTPUT: \>#3
        \end{tabbing}
                \begin{tabbing}
10.\ \=word\=word\=word\=word\=word\=word\=\kill
                \smallskip \setcounter{lineno}{0}
}
{           \end{tabbing}
            \par\nobreak\hrule\bigskip
            \rm
}

\newcommand{\nline}{\refstepcounter{lineno} %
\'\thelineno.\>%
}
\newcommand{\nnline}{\refstepcounter{lineno}%
\'\thelineno.\>\>%
}

\newcommand{\nnnline}{\refstepcounter{lineno}%
\'\thelineno.\>\>\>%
}

\newtheorem{theorem}{Theorem}[section]

\theoremstyle{definition}
\newtheorem{definition}[theorem]{Definition}

\theoremstyle{remark}
\newtheorem{remark}[theorem]{Remark}

\numberwithin{equation}{section}

\newcommand{\Z}{\mathbb{Z}}

\newcommand{\F}{\mathbb{F}}
\newcommand{\Q}{\mathbb{Q}}

\newcommand{\I}{\mathbb{I}}
\newcommand{\ie}{i.e., }

\newcommand{\ea}{{\em et al. }}
\renewcommand{\vec}[1]{\overline{\mathbf{#1}}}

\renewcommand{\labelitemi}{$-$}

\begin{document}

\title[Generalised Mersenne Numbers Revisited]{Generalised Mersenne Numbers Revisited}

\author{Robert~Granger}
\address{Claude Shannon Institute, Dublin City University, Ireland}
\curraddr{Claude Shannon Institute, UCD CASL, University College Dublin, Ireland}
\email{rgranger@computing.dcu.ie}
\thanks{The first author is supported by the Claude Shannon Institute, Science Foundation Ireland Grant No. 06/MI/006.}

\author{Andrew~Moss}
\address{Blekinge Institute of Technology, Sweden}
\email{awm@bth.se}

\subjclass[2010]{12Y05, 11T71, 11Y16}
%\MSC 12Yxx (Computational aspects of field theory and polynomials)
%\MSC 11T71 (Algebraic coding theory; cryptography) 
%\MSC 11Y16 (Algorithms; complexity)
\keywords{Prime fields, high-speed arithmetic, elliptic curve cryptography,
  generalised Mersenne numbers, cyclotomic primes, generalised repunit primes}

\begin{abstract}
Generalised Mersenne Numbers (GMNs) were defined by Solinas in 1999 and
feature in the NIST (FIPS 186-2) and SECG standards for use in elliptic curve
cryptography. Their form is such that modular reduction is extremely
efficient, thus making them an attractive choice for modular multiplication implementation.
However, the issue of residue multiplication efficiency seems to have been overlooked.
Asymptotically, using a cyclic rather than a linear convolution, residue
multiplication modulo a Mersenne number is twice as fast as integer
multiplication; this property does not hold for prime GMNs, unless they are
of Mersenne's form. In this work we exploit an alternative generalisation of Mersenne numbers
for which an analogue of the above property --- and hence the same efficiency ratio ---
holds, even at bitlengths for which schoolbook multiplication is optimal,
while also maintaining very efficient reduction. Moreover, our proposed
primes are abundant at any bitlength, whereas GMNs are extremely rare.
Our multiplication and reduction algorithms can also be easily
parallelised, making our arithmetic particularly suitable for hardware implementation. 
Furthermore, the field representation we propose also naturally protects against
side-channel attacks, including timing attacks, simple power analysis and
differential power analysis, which is essential in many cryptographic
scenarios, in constrast to GMNs.
\end{abstract}

\maketitle

\section{Introduction}
\label{sec:intro}

The problem of how to efficiently perform arithmetic in $\Z/N\Z$ is a very natural one,
with numerous applications in computational mathematics and number theory, such as primality 
proving~\cite{primality}, factoring~\cite{factoring}, and coding theory~\cite{coding}, for example. It is also of central 
importance to nearly all public-key cryptographic systems, including
the Digital Signature Algorithm~\cite{fips186-2}, RSA~\cite{RSA}, and elliptic curve cryptography (ECC)~\cite{ecc}.
As such, from both a theoretical and a practical perspective it is interesting and essential to have efficient 
algorithms for working in this ring, for either arbitrary or special moduli, with the application 
determining whether generality (essential for RSA for instance), or efficiency (desirable for ECC) takes precedence.

Two intimately related factors need consideration when approaching this problem. First, how should one represent residues? 
And second, how should one perform arithmetic on these representatives? A basic answer to the first
question is to use the canonical representation $\Z/N\Z = \{0,\ldots,N-1\}$. 
With regard to modular multiplication for example, an obvious answer to the second question
is to perform integer multiplication of residues, followed by reduction of the result modulo $N$,
in order to obtain a canonical representative once again.
Using this approach, the two components needed for efficient modular
arithmetic are clearly fast integer arithmetic, and fast modular reduction.

At bitlengths for which schoolbook multiplication is optimal, research on fast modular
multiplication has naturally tended to focus on reducing the cost of the reduction step. 
For arbitrary moduli, Montgomery's celebrated algorithm~\cite{mont}
enables reduction to be performed for approximately the cost of a residue by
residue multiplication. For the Mersenne numbers $M_k = 2^k - 1$, efficient modular multiplication consists of integer residue
multiplication to produce a $2k$-bit product $U \cdot 2^k + L$, with $U,L$ of
at most $k$-bits, followed by a single modular addition $U + L \bmod{M_k}$ to effect
the reduction, as is well known. 
In 1999 Solinas proposed an extension of this method to a larger class of integers: the Generalised
Mersenne Numbers (GMNs)~\cite{solinas}. As they are a superset, GMNs are more numerous than the Mersenne
numbers and hence contain more primes, yet incur little additional overhead in terms of performance~\cite{menezes}.
In 2000, NIST recommended ten fields for use in the ECDSA: five binary fields and five prime
fields, and due to their performance characteristics the latter of these are all GMNs~\cite{fips186-2},
which range from $192$ to $521$ bits in size. The Standards for Efficient Cryptography
Group also recommended the same five prime fields in 2010~\cite{sec2}.

For the GMNs recommended by NIST, there is no interplay between the
residue multiplication and reduction algorithms, each step being treated separately
with respect to optimisation.
On the other hand, at asymptotic bitlengths the form of the modulus may be effectively
exploited to speed up the residue multiplication step. For the Mersenne
numbers $M_k$ in particular, modular multiplication can be performed for any $k$ using a cyclic convolution effected by an
irrational-base discrete weighted transform (IBDWT)~\cite[\S6]{crandallfagin} (see
also~\cite[\S9.5.2-9.5.4]{crandallpom} for an excellent overview of discrete Fourier
transform-based multiplication methods, convolution theory and IBDWTs).  
As such, multiplication modulo Mersenne numbers is approximately twice as fast as
multiplication of integers of the same bitlength, for which a linear
convolution is required, as each multiplicand must be padded with $k$ zeros
before a cyclic convolution of length $2k$ can be performed.
For Montgomery multiplication at asymptotic bitlengths, the reduction step can be made $25\%$ cheaper, 
again by using a cyclic rather than a linear convolution for one of the
required multiplications~\cite{phatak}. However, since the multiplication step
is oblivious to the form of the modulus, it seems unlikely to possess the
same efficiency benefits that the Mersenne numbers enjoy. These considerations raise the
natural question of whether there exists a similar residue multiplication speed-up at bitlengths for which
schoolbook multiplication is optimal? Certainly for the modulus $N = 2^k$, such a speed-up can be achieved, since
the upper half words of the product can simply be ignored. However, this
modulus is unfortunately not at all useful for ECC.

In this work we answer the above question affirmatively, using an alternative generalisation
of Mersenne numbers, which has several desirable features:
\vspace{-2mm}
\begin{list}{\labelitemi}{\leftmargin=1.2em}
\item {\bf Simple.} Our proposed family is arguably a far more natural generalisation of Mersenne numbers than Solinas', 
and gives rise to beautiful multiplication and reduction algorithms.
\item {\bf Abundant.} Our primes are significantly more numerous than the set of prime GMNs and
are abundant for all tested bitlengths; indeed their number can be estimated using Bateman and
Horn's quantitative version~\cite{bateman} of Schinzel and Sierpi\'{n}ski's ``Hypothesis H''~\cite{schinzel}. 
\item {\bf Fast multiplication.} Our residue multiplication is nearly twice as fast as multiplication of integer residues.
\item {\bf Fast reduction.} Our reduction has linear complexity and is
  particularly efficient for specialised parameters, although such
  specialisation comes at the cost of reducing the number of primes available. 
\item {\bf Parallelisable.} Both multiplication and reduction can be easily parallelised, 
making our arithmetic particularly suitable for hardware implementation. 
\item {\bf Side-channel secure.} Our representation
naturally protects against well-known side-channel attacks on ECC (see~\cite[ch. IV]{ecc2} for an
overview), in contrast to the NIST GMNs, see~\cite{sakai} and~\cite[\S3.2]{danrandom}. 
This includes timing attacks~\cite{kocher1,walter2}, simple power
analysis~\cite{sakai} and differential power analysis~\cite{kocher2}. 
\end{list}
\vspace{+1mm}

This article provides an introductory (and comprehensive) theoretical
framework for the use of our proposed moduli. It thus serves as a foundation for a new
approach to the secure and efficient implementation of prime fields for ECC, both in software and in hardware.
At a high level, our proposal relies on the combination of a remarkable algebraic identity
used by Nogami, Saito, and Morikawa in the context of extension fields~\cite{nogami}, together with the 
residue representation and optimisation of the reduction method proposed by Chung and Hasan~\cite{CH3}, which models 
suitable prime fields as the quotient of an integer lattice by a particular equivalence relation.
To verify the validity of our approach, we also provide a proof-of-concept 
implementation that is already competitive with the current fastest modular
multiplication algorithms at contemporary ECC security levels~\cite{25519,speed,hisil,longa,GLS,bernstein2}.

The sequel is organised as follows. In \S\ref{sec:definitions} we present some
definitions and recall related work. In \S\ref{sec:MRCPfieldrep} we
describe the basis of our arithmetic, then in
\S\ref{sec:multiplication}-\ref{sec:representation} we present details of our
residue multiplication, reduction and representation respectively.
In \S\ref{sec:stability} we show how to ensure I/O stability for modular
multiplication, then in \S\ref{sec:fullmul} we put everything together into a
full modular multiplication algorithm.
We then address other arithmetic operations and give a brief treatment of side-channel
secure ECC in \S\ref{sec:otherops}, and in \S\ref{sec:paramgen} show how to 
generate suitable parameters. In \S\ref{sec:results} we present our
implementation results and finally, in \S\ref{sec:conclusion} we draw some conclusions.

%-------------------------------------------------------------------------------------------------------------------------
%-------------------------------------------------------------------------------------------------------------------------

\section{Definitions and Related Work}
\label{sec:definitions}

In this section we introduce the cyclotomic primes and provide a summary of
related work. We begin with the following definition.

\begin{definition}
For $n \geq 1$ let $\zeta_{n}$ be a primitive $n$-th root of unity.
The $n$-th cyclotomic polynomial is defined by
\begin{equation}
\nonumber \Phi_n(x) = \prod_{(k,n) = 1} (x - \zeta_{n}^{k}) = \prod_{d|n} (1-x^{n/d})^{\mu(d)},
\end{equation}
where $\mu$ is the M\"{o}bius function.
\end{definition}

Two basic properties of the cyclotomic polynomials are that they
have integer coefficients, and are irreducible over $\Z$.
These two properties ensure that the evaluation of a cyclotomic polynomial
at an integer argument will also be an integer, and that this integer will not inherit a
factorisation from one in $\Z[x]$. One can therefore ask whether or not these
polynomials ever assume prime values at integer arguments, which leads to our next definition.

\begin{definition}
For $n \geq 1$ and $t \in \Z$, if $p = \Phi_n(t)$ is prime, we call $p$ an $n$-th cyclotomic prime,
or simply a cyclotomic prime.
\end{definition}

Note that for all primes $p$, we have $p = \Phi_1( p + 1 ) = \Phi_2(p-1)$, and so trivially all
primes are cyclotomic primes. These instances are also trivial in the context of the algorithms we present for performing
arithmetic modulo these primes, since in both cases
the cyclotomic polynomials are linear and our algorithms reduce to ordinary Montgomery arithmetic.
Hence for the remainder of the article we assume $n \geq 3$.

In addition to being prime-evaluations of cyclotomic polynomials, note that for a cyclotomic prime $p=\Phi_n(t)$,
the field $\F_p$ can be modelled as the quotient of the ring of integers of
the $n$-th cyclotomic field $\Q(\zeta_n)$, by the prime ideal $\pi = \langle p, \zeta_n - t \rangle$.
This is precisely how one would represent $\F_p$ when applying the Special
Number Field Sieve to solve discrete logarithms in $\F_p$, for example~\cite{LenstraBook}.
Hence our nomenclature for these primes seems apt. This interpretation of $\F_p$ for $p$ a cyclotomic prime is implicit within
the arithmetic we develop here, albeit only insofar as it provides a theoretical context for it; 
this perspective offers no obvious insight into how to perform arithmetic efficiently and the algorithms we develop make no
use of it at all. Similarly, the method of Chung and Hasan~\cite{CH3} upon which our residue representation is based 
can be seen as arising in exactly the same way for the much larger set of primes they consider,
with the field modelled as a quotient of the ring of integers of a suitable number field by a degree one prime ideal,
just as for the cyclotomic primes.  

\subsection{Low redundancy Cyclotomic Primes}
\label{subsec:lrcp}

The goal of the present work is to provide efficient algorithms for performing
$\F_p$ arithmetic, for $p = \Phi_n(t)$ a cyclotomic prime. As will become clear from our
exposition, in order to exploit the available cyclic structure --- for both multiplication and reduction --- we do
not use the field $\Z/\Phi_{n}(t)\Z$, but instead embed into the slightly
larger ring $\Z/(t^{n} - 1)\Z$ if $n$ is odd, and $\Z/(t^{n/2} + 1)\Z$ if $n$
is even. In each case, using the larger ring potentially introduces an
expansion factor $e(n)$ into the residue representation. One can alternatively
view this in terms of a redundancy measure $r(n)$, where $r = e-1$. Since
using a larger ring for arithmetic will potentially be slower, we now identify
three families of cyclotomic polynomials for which the above embeddings have low redundancy.

For $n$ even, there is a family of cases for which the above embedding
does not introduce any redundancy, namely for $n = 2^k$, since
$\Phi_{2^k}(t) = t^{2^{k-1}} + 1 = t^{2^k/2} + 1$, and hence $e=1$ and $r=0$. 
When $t=2$ these are of course the Fermat numbers, and for general $t$ these 
integers are known as Generalised Fermat Numbers (GFNs). It is expected that 
for each $k$ there are infinitely many $t$ for which $t^{2^k} + 1$ is prime~\cite[\S3]{dubner}.

If $n = 2p$ for $p$ prime, then $\Phi_{2p}(t) = t^{p-1} - t^{p-2} + \cdots +
t - 1 = (t^{p} + 1)/(t+1)$ and in this case $e = p/(p-1)$ and $r = 1/(p-1)$. 
The primality of these numbers was studied in~\cite{dubner2}, and while
they apparently do not have a designation in the literature, one can see that
by substituting $t$ with $-t$ in the third family below produces this one.
For general even $n$ we have $e = n/2\phi(n)$ and $r = (n-2\phi(n))/2\phi(n)$,
with $\phi(\cdot)$ Euler's totient function, which is the degree of $\Phi_n(x)$.
Hence amongst those even $n$ which are not a power of $2$, this family
produces the successive local minima of $r$.

For odd $n$, we have $e = n/\phi(n)$ and $r = (n - \phi(n))/\phi(n)$. The successive
local minima of $r$ occur at $n = p$ for $p$ prime, in which case 
$\Phi_{p}(t) = t^{p-1} + t^{p-2} + \cdots + t + 1 = (t^p - 1)/(t-1)$, also with $r = 1/(p-1)$.
When $t=2$ these are of course the Mersenne numbers, and in analogy with the
case of Fermat numbers, it would be natural to refer to these integers for general
$t$ as Generalised Mersenne Numbers, particularly as one can show they share the
aforementioned asymptotic efficiency properties of the Mersenne numbers, while
Solinas' GMNs do not, unless they are of Mersenne's form.
However, this family of numbers is known in the literature as generalised 
repunits~\cite{williams,snyder,dubner3}, since their base-$t$ expansion
consists entirely of $1$'s. 
Therefore for the sake of uniform nomenclature, we use the following definition.

\begin{definition}
For $m+1$ an odd prime let
\[
p = \Phi_{m+1}(t) = t^m + t^{m-1} + \cdots + t + 1.
\]
We call such an integer a {\em Generalised Repunit (GR)}; when $p$ is prime we
call it a {\em Generalised Repunit Prime (GRP)}.
\end{definition}

We have developed modular multiplication algorithms for both GRPs and GFNs.
In terms of efficiency, for GRPs and GFNs of the
same bitlength the respective multiplication algorithms require exactly
the same number of word-by-word multiplications. 
Also, our reduction algorithms for both GRPs and GFNs are virtually identical.
However, the multiplication algorithm for GFNs is far less elegant, is not perfectly parallelisable and
contains more additions. Furthermore, for a given bitlength there are fewer efficient GFN primes than there are
GRPs --- as the bitlength of GFNs doubles as $k$ is incremented --- and the I/O
stability analysis for multiplication modulo a GRP is far simpler. Therefore in this exposition we focus on algorithms for performing
arithmetic modulo GRPs and their analysis only. Note that the studies of GRPs~\cite{dubner3,williams} consider only very small $t$ and
large $m$, whereas we will be interested in $t$ approximately the word base of
the target architecture, and $m$ the number of words in the prime whose field arithmetic
we are to implement. Hence one expects (and finds) there to be very many GRPs
for any given relevant bitlength, see \S\ref{sec:paramgen}.

\subsection{Related work}\label{subsec:previouswork}

In the context of extension fields, let $m+1$ be prime and let $p$ be a
primitive root modulo $m+1$. Then $\F_{p^m} =  \F_p[x]/(\Phi_{m+1}(x)\F_p[x])$.
In the binary case, \ie $p=2$, several authors have proposed the use of
this polynomial --- also known as the all-one polynomial (AOP) --- 
to obtain efficient multiplication algorithms~\cite{itoh_tsujii,wolf,blake,silverman}. 
All of these rely on the observation that the field
$\F_{2}[x]/(\Phi_{m+1}(x)\F_{2}[x])$ embeds into the ring
$\F_{2}[x]/((x^{m+1} + 1)\F_{2}[x])$ --- referred to by Silverman~\cite{silverman} as the
``ghost bit'' basis --- which possesses a particularly
nice cyclic structure, but introduces some redundancy. Similarly, this idea applies to any
cyclotomic polynomial, and several authors have investigated this
strategy, embedding suitably defined extension fields into
the ring $\F_{2}[x]/((x^{n} + 1)\F_{2}[x])$~\cite{drolet,geiselmann,hasan}. 

For odd characteristic extension fields, Silverman noted that the ``ghost bit'' basis for
$p=2$ extends easily to larger $p$~\cite{silverman}, while Kwon~\ea have explored this idea further~\cite{kwon}.
Central to our application is the work of Nogami, Saito and Morikawa~\cite{nogami},
who used the AOP to obtain a very fast multiplication
algorithm, see~\S\ref{sec:multiplication}.
The use of cyclotomic polynomials in extension field arithmetic is therefore well studied.
In the context of prime fields however, the present work appears to be the first to
transfer ideas for cyclotomic polynomials from the domain of extension field arithmetic to prime
field arithmetic, at least for the relatively small bitlengths for which
schoolbook multiplication is optimal.

With regard to the embedding of a prime field into a larger integer ring, 
the idea of operand scaling was introduced by Walter
in order to obtain a desired representation in the higher-order bits~\cite{walter},
which aids in the estimation of the quotient when using Barrett reduction~\cite{barrett}.
Similarly, Ozturk~\ea proposed using fields with characteristics dividing
integers of the form $2^k \pm 1$, with particular application to ECC~\cite{ozturk}.
As stated in the introduction, there are numerous very efficient prime field
ECC implementations~\cite{25519,speed,hisil,longa,bernstein2}.
While the moduli used in these instances permit fast reduction algorithms, and
the implementations are highly optimised, it would appear that none of
them permit the same residue multiplication speed-up that we present here,
which is one of the central distinguishing features of the present work.

%-------------------------------------------------------------------------------------------------------------------------
%-------------------------------------------------------------------------------------------------------------------------

\section{GRP Field Representation}
\label{sec:MRCPfieldrep}

In this section we present a sequence of representations of $\F_p$, with $p$
a GRP, the final one being the target representation which we use for our arithmetic.
We recall the mathematical framework of Chung-Hasan arithmetic, in both the
general setting and as specialised to GRPs, focusing here on the underlying
theory, deferring explicit algorithms for residue multiplication, reduction and representation until 
\S\ref{sec:multiplication}-\ref{sec:representation}.

\subsection{Chung-Hasan arithmetic}
\label{subsec:CHarithmetic}

We now describe the ideas behind Chung-Hasan arithmetic~\cite{CH1,CH2,CH3}. 
The arithmetic was developed
for a class of integers they term low-weight polynomial form integers (LWPFIs),
whose definition we now recall.

\begin{definition}
An integer $p$ is a low-weight polynomial form integer (LWPFI), if it can be represented
by a monic polynomial $f(t) = t^n + f_{n-1} t^{n-1} + \cdots + f_1 t + f_0$, where $t$ is a positive
integer and $|f_i| \leq \xi$ for some small positive integer $\xi < t$.
\end{definition}

Note that if for a given LWPFI each $f_i \in \{\pm1,0\}$ and $t = 2^k$, then it is a GMN, as defined by Solinas~\cite{solinas}.
The key idea of Chung and Hasan is to perform arithmetic
modulo $p$ using representatives from the polynomial ring $\Z[T]/(f(T)\Z[T])$.
To do so, one uses the natural embedding
$\psi: \F_p \hookrightarrow \Z[T]/(f(T)\Z[T])$ obtained
by taking the base $t$ expansion of an element of $\F_p$ in the
canonical representation $\F_p = \{0,\ldots,p-1\}$, and substituting $T$ for $t$.
To compute $\psi^{-1}$ one simply makes the inverse substitution and evaluates
the expression modulo $p$.

The reason for using this ring is straightforward: since $\psi^{-1}$ is a homomorphism,
when one computes $z(T) = x(T) \cdot y(T)$ in $\Z[T]$, reducing the result modulo $f(T)$
to give $w(T)$ does not change the element of $\F_p$ represented by $z(T)$, \ie if
$z(T) \equiv w(T) \pmod {f(T)}$, then $z(t) \equiv w(t) \pmod{p}$, since $p=f(t)$.
Furthermore, since $f(T)$ has very small coefficients, $w(T)$
can be computed from $z(T)$ using only additions and subtractions.
Hence given the degree $2(n-1)$ product of two degree $n-1$ polynomials in $\Z[T]$,
its degree $n-1$ representation in $\Z[T]/(f(T)\Z[T])$ can be computed very efficiently.
Note that for non-low-weight polynomials this would no longer be the case.

The only problem with this approach is that when computing
$z(T)$ as above, the coefficients of $z(T)$, and hence $w(T)$,
will be approximately twice the size of the inputs' coefficients,
and if further operations are performed the representatives will continue
to expand. Since for I/O stability one requires that the coefficients
be approximately the size of $t$ after each modular multiplication or squaring,
one must somehow reduce the coefficients of $w(T)$ to obtain a
standard, or reduced representative, while ensuring that $\psi^{-1}(w(T))$ remains unchanged.

Chung and Hasan refer to this issue as the
{\em coefficient reduction problem} (CRP), and developed three solutions in their
series of papers on LWPFI arithmetic~\cite{CH1,CH2,CH3}. Each of these solutions
is based on an underlying lattice, although this was only made explicit in~\cite{CH3}.
Since the lattice interpretation is the most elegant and simplifies the exposition, 
in the sequel we opt to develop the necessary theory for GRP arithmetic in this setting.

\subsection{Chung-Hasan representation for GRPs}
\label{subsec:overview}

Let $p = \Phi_{m+1}(t)$ be a GRP. Our goal is to develop arithmetic for
$\F_p$, and we begin with the canonical representation $\F_p = \Z/\Phi_{m+1}(t)\Z$. 
As stated in \S\ref{subsec:lrcp}, the first map in our chain of representations takes the canonical ring and embeds it into
$\Z/(t^{m+1}-1)\Z$, for which the identity map suffices. To map back, one reduces a representative
modulo $p$. We then apply the Chung-Hasan transformation of \S\ref{subsec:CHarithmetic}, which embeds the second ring into
$\Z[T]/(T^{m+1}-1)\Z[T]$, by taking the base $t$ expansion of a canonical
residue representative in $\Z/(t^{m+1}-1)\Z$, and substituting $T$ for $t$. We call this map $\psi$. 
To compute $\psi^{-1}$ one simply makes the inverse substitution and evaluates the expression modulo $t^{m+1} - 1$.

Note that the codomain of $\psi$ may be regarded as an $(m+1)$-dimensional vector space over $\Z$, 
equipped with the natural basis $\{T^m,\ldots,T,1\}$. In particular, for $x(T) \in \Z[T]/(T^{m+1}-1)\Z[T]$, where
\[
x(T) = x_mT^m + \ldots +x_1T + x_0,
\]
one can consider $x(T)$ to be a vector $\vec{x} = [x_m,\ldots,x_0] \in \Z^{m+1}$.
Since $\Z^{m+1}$ has elements whose components are naturally unbounded, 
for each $x \in \Z/(t^{m+1}-1)\Z$ there are infinitely many elements of
$\Z^{m+1}$ that map via $\psi^{-1}$ to $x$. Therefore in order to obtain a useful isomorphism directly between
$\Z/(t^{m+1}-1)\Z$ and $\Z^{m+1}$, we identify two elements of $\Z^{m+1}$
whenever they map via $\psi^{-1}$ to the same element of $\Z/(t^{m+1}-1)\Z$, \ie
\begin{equation}\label{eq:lattice_mod}
\vec{x} \sim \vec{y} \hspace{2mm} \Longleftrightarrow \hspace{2mm}
\psi^{-1}(\vec{x}) \equiv \psi^{-1}(\vec{y}) \pmod {t^{m+1} - 1},
\end{equation}
and take the image of $\psi$ to be the quotient of $\Z^{m+1}$ by this equivalence relation. 
Pictorially, we thus have: 
\[
\F_p \subset \Z/(t^{m+1} - 1)\Z \cong \Z^{m+1}/\sim
\]

As mentioned in \S\ref{subsec:CHarithmetic}, for each coset in $\Z^{m+1}/\sim$, 
we should like to use a minimal, or in some sense `small' representative,
in order to facilitate efficient arithmetic after a multiplication or a squaring, for example. 
Since we know that the base-$t$ expansion of every $x \in \Z/(t^{m+1}-1)\Z$ 
gives one such representative for each coset in $\Z^{m+1}/\sim$, for a reduction 
algorithm we just need to be able to find it, or at least one whose components
are of approximately the same size. Chung and Hasan related finding such
`nice' or reduced coset representatives to solving a computational
problem in an underlying lattice, which we now recall.

\subsection{Lattice interpretation}
\label{subsec:lattice}

Given an input vector $\vec{z}$, which is the output of a multiplication or a squaring, a coefficient reduction algorithm should output
a vector $\vec{w}$ such that $\vec{w} \sim \vec{z}$, in the sense of~(\ref{eq:lattice_mod}),
whose components are approximately the same size as $t$. As observed in~\cite{CH3},
the equivalence relation~(\ref{eq:lattice_mod}) is captured
by an underlying lattice, and finding $\vec{w}$ is tantamount
to solving an instance of the {\em closest vector problem} (CVP) in this lattice.
To see why this is, we first fix some notation as in~\cite{CH3}.

Let $\vec{u}$ and $\vec{v}$ be vectors in $\Z^{m+1}$ such that the following condition
is satisfied:
\begin{equation}\label{relation}
\nonumber [t^{m}, \ldots, t, 1] \cdot \vec{u}^T \equiv [t^{m}, \ldots, t, 1] \cdot \vec{v}^T \pmod {t^{m+1}-1}
\end{equation}
Then we say that {\em $\vec{u}$ is congruent to $\vec{v}$ modulo $t^{m+1}-1$} and write this as
$\vec{u} \cong_{t^{m+1}-1} \vec{v}$. Note that this is exactly the same as saying
$\psi^{-1}(\vec{u}) \equiv \psi^{-1}(\vec{v}) \pmod {t^{m+1}-1}$, and so
$\vec{u} \sim \vec{v} \hspace{0mm} \Longleftrightarrow \hspace{0mm} \vec{u} \cong_{t^{m+1}-1} \vec{v}$.

Similarly, but abusing notation slightly,
for any integer $b \neq t^{m+1}-1$ (where $b$ is typically a power of the word base of the target architecture),
we write $\vec{u} \cong_{b} v$ for some
integer $v$ satisfying $[t^{m}, \ldots, t, 1] \cdot \vec{u}^T \equiv v \pmod {b}$,
and say {\em $\vec{u}$ is congruent to $v$ modulo $b$}, in this case.
We reserve the use of `$\equiv$' to express a component-wise congruence relation, \ie
$\vec{u} \equiv \vec{v} \pmod {b}$. 
Finally, we denote by $\vec{u} \bmod b$ the component-wise modular reduction of $\vec{u}$ by $b$.

The lattice underlying the equivalence relation~(\ref{eq:lattice_mod}) can now enter the frame.
Let $\mathbf{V} = \{\vec{v}_0,\ldots,\vec{v}_{m}\}$ be a set of $m+1$ linearly independent vectors in $\Z^{m+1}$
such that $\vec{v}_i \cong_{t^{m+1}-1} \vec{0}$, the all zero vector, for $i=0,\ldots,m$. Then the set of all integer
combinations of elements of $\mathbf{V}$ forms an integral lattice, $\mathcal{L}(\mathbf{V})$, with the
property that for all $\vec{z} \in \Z^{m+1}$, and all $\vec{u} \in \mathcal{L}$, we have
\begin{equation}\label{eq:no_sum}
\vec{z} + \vec{u} \cong_{t^{m+1}-1} \vec{z}
\end{equation}
In particular, the equivalence relation~(\ref{eq:lattice_mod}) is captured by the
lattice $\mathcal{L}$, in the sense that
\begin{equation}\label{eq:new_equiv}
\nonumber \vec{x} \cong_{t^{m+1}-1} \vec{y} \hspace{2mm} \Longleftrightarrow \hspace{2mm}
\vec{x} - \vec{y} \in \mathcal{L}
\end{equation}
Therefore if one selects basis vectors for $\mathcal{L}$ that have infinity-norm approximately $t$,
then for a given $\vec{z} \in \Z^{m+1}$, finding the closest vector $\vec{u} \in \mathcal{L}$ to $\vec{z}$ 
(with respect to the $L_{\infty}$-norm), means the vector $\vec{w} = \vec{z} - \vec{u}$
is in the fundamental domain of $\mathcal{L}$, and so
has components of the desired size. Furthermore, since $\vec{w} = \vec{z} - \vec{u}$,
by~(\ref{eq:no_sum}) we have 
\[
\vec{w} \cong_{t^{m+1}-1} \vec{z},
\]
and hence solving the CVP in this lattice solves the CRP. In general solving the CVP is NP-hard,
but since we can exhibit a good (near-othogonal) lattice basis for LWPFIs, and
an excellent lattice basis for GRPs, solving it is straightforward in our case.

\subsection{Lattice basis and simple reduction}
\label{subsec:latticereduction}

For GRPs, we use the following basis for $\mathcal{L}$:
\begin{equation}\label{eq:latticebasis}
\left[ \begin{array}{cccccc}
1  & 0 & \cdots & 0 & 0 & -t \\
-t & 1 & \cdots & 0 & 0 & 0   \\
0 & -t & \cdots & 0 & 0 & 0   \\
\vdots & \vdots & \ddots & \vdots &\vdots & \vdots \\
0 & 0 & \cdots & -t & 1 & 0 \\
0 & 0 & \cdots & 0 & -t & 1
\end{array} \right]
\end{equation}

Observe that the infinity-norm of each basis vector is $t$, so elements in the fundamental 
domain will have components of the desired size, and that each basis vector is
orthogonal to all others except the two adjacent vectors (considered cyclically).
In order to perform a simple reduction that reduces the size of components by approximately $\log_2t$ bits,
write each component of $\vec{z}$ in base $t$: $z_i = z_{i,1}t + z_{i,0}$. If we define $\vec{w}^T$ to be: 
\[
\left[ \begin{array}{c}
%    z_{m,1} t + z_{m,0}\\
%    z_{m-1,1} t + z_{m-1,0}\\
%    \vdots\\
%    \vdots\\
%    z_{1,1} t + z_{1,0}\\
%    z_{0,1} t + z_{0,0}
%    \end{array} \right]
   z_{m}\\
   z_{m-1}\\
   \vdots\\
   \vdots\\
   z_{1}\\
   z_{0}
   \end{array} \right]
+
\left[ \begin{array}{cccccc}
1  & 0 & \cdots & 0 & 0 & -t \\
-t & 1 & \cdots & 0 & 0 & 0   \\
0 & -t & \cdots & 0 & 0 & 0   \\
\vdots & \vdots & \ddots & \vdots &\vdots & \vdots \\
0 & 0 & \cdots & -t & 1 & 0 \\
0 & 0 & \cdots & 0 & -t & 1
\end{array} \right]
\left[ \begin{array}{c}
    z_{m-1,1}\\
    z_{m-2,1}\\
    \vdots\\
    \vdots\\
    z_{0,1}\\
    z_{m,1}\\
    \end{array} \right],
\]
then $\vec{w} \cong_{t^{m+1}-1} \vec{z}$ and each $|w_i| \approx |z_i| / t$, assuming $|z_i| > t^2$. This was the method of reduction 
described in~\cite{CH1}, which requires integer division. The idea described in~\cite{CH2} was
based on an analogue of Barrett reduction~\cite{barrett}. The method we shall use, from~\cite{CH3}, is based
on Montgomery reduction~\cite{mont} and for $t$ not a power of $2$ is the most
efficient of the three Chung-Hasan methods.

\subsection{Montgomery lattice-basis reduction}\label{subsec:montlattice}

In ordinary Montgomery reduction~\cite{mont}, one has an integer $0 \leq Z < pR$ which
is to be reduced modulo $p$, an odd prime, where here $R$ is the smallest
power of the word base $b$ larger than $p$. 
The central idea is to add a multiple of $p$ to $Z$ such that the result
is divisible by $R$. Upon dividing by $R$, which is a simple right shift of words,
the result is congruent to $ZR^{-1} \pmod {p}$, and importantly is less than $2p$.

In the context of GRPs, let $R = b^q $ be the smallest power of $b$ greater than $t$.
The input to the reduction algorithm is a vector $\vec{z} \in \Z^{m+1}$ for which each component
is approximately $R^2$. The natural analogue of Montgomery reduction is to
add to $\vec{z}$ a vector $\vec{u} \in \mathcal{L}$ whose components
are also bounded by $R^2$, such that $\vec{z} + \vec{u} \equiv [0,\ldots,0] \pmod{R}$.
Then upon the division of each component by $R$, the result will be a vector $\vec{w}$ which
satisfies
\[
\vec{w} \cong_{t^{m+1}-1}  (\vec{z} + \vec{u}) \cdot R^{-1} \cong_{t^{m+1}-1} \vec{z} \cdot R^{-1},
\]
and which has components of the desired size.
While this introduces an $R^{-1}$ term into the congruence, as with Montgomery arithmetic,
one circumvents this simply by altering the original coset representation of $\Z/(t^{m+1}-1)\Z$,
via the map $x \mapsto xR \pmod {t^{m+1}-1}$, which is bijective since $\gcd(t^{m+1}-1,R) = 1$,
assuming $t$ is even, see~\S\ref{sec:reduction}.
How then does one find a suitable lattice point $\vec{u}$? For this one use the 
lattice basis~(\ref{eq:latticebasis}), which from here on in we call $L$.
Proposition 3 of~\cite{CH3} proves that $\det L = 1-t^{m+1}$, and so $\gcd(\det L, R)=1$. 
One can therefore compute
\begin{eqnarray}
\label{basic_mont1} \vec{u}^{T} &\stackrel{{\rm def}}{=} & -L^{-1} \cdot \vec{z}^{T} \pmod{R},\\
\label{basic_mont2} \vec{w}^T &\stackrel{{\rm def}}{=}& (\vec{z}^T + L \cdot \vec{u}^T)/R,
\end{eqnarray}
giving $\vec{w}$ with the required properties.
Observe that the form of these two operations is identical to Montgomery reduction, the only 
difference being that integer multiplication is replaced by matrix by vector multiplication.
It is easy to see that this is what one requires, since for any $\vec{u} \in \Z^{m+1}$, we have
$L \cdot \vec{u}^T \in \mathcal{L}$, and so 
\[
\vec{z}^T + L \cdot \vec{u}^T \cong_{t^{m+1}-1} \vec{z}^T.
\]
Furthermore, modulo $R$ we have
\begin{equation}
\nonumber \vec{z}^T + L \cdot \vec{u}^T = \vec{z}^T + L \cdot (-L^{-1} \cdot \vec{z}^T \bmod R) \equiv
[0,\ldots,0]^T,
\end{equation}
ensuring the division of each component by $R$ is exact. Hence $\vec{w} \cong_{t^{m+1}-1} \vec{z} \cdot R^{-1}$, 
as claimed.

In~\cite{CH3}, an algorithm was given for computing $\vec{u}$ and $\vec{w}$ in~(\ref{basic_mont1})
and~(\ref{basic_mont2}) respectively, for an arbitrary LWPFI $f(t)$.
The number of word-by-word multiply instructions in the algorithm --- which is
the dominant cost --- is $\approx n q^2$, 
where $n$ is the degree of $f(t)$, and $R = b^q$. In comparison, for ordinary Montgomery reduction 
modulo an integer of equivalent size this number is $n^2 q^2$, making the former approach potentially 
very attractive. For our choice of primes --- the GRPs --- our specialisation of this algorithm
is extremely efficient, as we show in \S\ref{sec:reduction}.

\subsection{High level view of Chung Hasan-arithmetic}
\label{subsec:highlevel}

For extension fields, there exists a natural separation between the polynomial
arithmetic of the extension, and the prime subfield arithmetic, which makes respective
optimisation considerations for each almost orthogonal. On the other hand, if
for an LWPFI one naively attempts to use efficient techniques that are valid for extension fields, 
then one encounters an inherent obstruction, namely that there is no such separation between the
polynomial arithmetic and the coefficient arithmetic, which leads to
coefficient expansion upon performing arithmetic operations. Chung-Hasan
arithmetic can be viewed as a tool to overcome this obstruction, since it
provides an efficient solution to the coefficent reduction problem. In
practice therefore any efficient techniques for extension field
arithmetic can be ported to prime fields, whenever the prime is an LWPFI,
which is precisely what we do in \S\ref{sec:multiplication}.

%-------------------------------------------------------------------------------------------------------
%-------------------------------------------------------------------------------------------------------

\section{GRP Multiplication}
\label{sec:multiplication}

In this section we detail algorithms for performing multiplication of GRP
residue representatives. While for the reduction and residue representation we consider elements to be in $\Z^{m+1}$, the
multiplication algorithm arises from the arithmetic of the polynomial ring $\Z[T]/(T^{m+1}-1)\Z[T]$,
and so here we use this ring to derive the multiplication formulae.

\subsection{Ordinary multiplication formulae}
\label{subsec:mulordinary}

Let $\mathcal{R} = \Z[T]/(T^{m+1}-1)\Z[T]$, and let $\vec{x} = [x_m,\ldots,x_0]$ and $\vec{y}=[y_m,\ldots,y_0]$
be elements in $\mathcal{R}$. Then in $\mathcal{R}$ the product $\vec{x} \cdot \vec{y}$ is equal to
$[z_m,\ldots,z_0]$, where
\begin{equation}\label{eq:convolution}
z_i = \sum_{j=0}^{m} x_{\langle j \rangle} y_{\langle i-j \rangle},
\end{equation}
where the subscript $\langle i \rangle$ denotes $i \pmod{m+1}$.  
This follows from the trivial property $T^{m+1} \equiv 1 \pmod{T^{m+1}-1}$,
and that for $\vec{x} = \sum_{i=0}^{m} x_i T^i$ and $\vec{y} = \sum_{j=0}^{m} y_j T^j$, we have:
\begin{eqnarray}
\nonumber \vec{x} \cdot \vec{y} &=& \sum_{i=0}^{m} x_i \cdot (T^i \cdot
\vec{y}) = \sum_{i=0}^{m} x_i \cdot \Big(\sum_{j=0}^{m} y_j T^{i+j}\Big)\\
\nonumber &=& \sum_{i=0}^{m} x_i \cdot \Big(\sum_{j=0}^{m} y_{\langle j-i
\rangle } T^{j}\Big) = \sum_{j=0}^{m} \Big(\sum_{i=0}^{m} x_i \cdot y_{
\langle j-i \rangle }\Big) T^{j}.
\end{eqnarray}
This is of course just the cyclic convolution of $\vec{x}$ and $\vec{y}$.

\subsection{Multiplication formulae of Nogami \ea}\label{subsec:mulnogami}

Nogami, Saito and Morikawa proposed the use of all-one polynomials (AOPs) to
define extensions of prime fields~\cite{nogami}. In this section we
will first describe their algorithm in this context, and then show how it fits into the
framework developed in~\S\ref{sec:MRCPfieldrep}. 

Let $\F_p$ be a prime field and let $f(\omega) = \omega^m + \omega^{m-1} + \cdots + \omega + 1$ be irreducible over
$\F_p$, \ie $m+1$ is prime and $p$ is a primitive root modulo $m+1$.
Then $\F_{p^m} = \F_p[\omega]/(f(\omega)\F_p[\omega])$.
Using the polynomial basis $\{\omega^m,\omega^{m-1},\ldots,\omega\}$ ---
rather than the more conventional $\{\omega^{m-1},\ldots,\omega,1\}$ ---
elements of $\F_{p^m}$ are
represented as vectors of length $m$ over $\F_p$:
\[
\vec{x} = [x_m,\ldots,x_1] = x_m\omega^m + x_{m-1}\omega^{m-1} + \cdots + x_1\omega.
\]
Let $\vec{x} = [x_m,\ldots,x_1]$ and $\vec{y} = [y_m,\ldots,y_1]$
be two elements to be multiplied. For $0 \leq i \leq m$, let
\begin{equation}\label{nogami_mul}
q_i = \sum_{j=1}^{m/2} (x_{\langle \frac{i}{2} + j\rangle} - x_{\langle \frac{i}{2} - j\rangle})
(y_{\langle \frac{i}{2} + j\rangle} - y_{\langle \frac{i}{2} - j\rangle}),
\end{equation}
where the subscript $\langle i \rangle$ here, as in \S\ref{subsec:mulordinary}, denotes $i \pmod{m+1}$.
One then has:
\begin{equation}\label{eq:nogami_state}
\vec{z} = \vec{x} \cdot \vec{y} = \sum_{i=1}^{m} z_i \omega^i, \ \ \text{with} \ z_i = q_0 - q_i.
\end{equation}
Nogami \ea refer to these coefficient formulae as the
{\em cyclic vector multiplication algorithm} (CVMA) formulae. The CVMA formulae are remarkable, 
since the number of $\F_p$ multiplications is reduced relative to the schoolbook method
from $m^2$ to $m(m+1)/2$, but at the cost of increasing the number of $\F_p$ additions from $m^2 - 1$ to
$3m(m-1)/2 - 1$. As alluded to in \S\ref{subsec:highlevel}, a basic insight of
the present work is the observation that one may apply the expressions in~(\ref{nogami_mul})
to GRP multiplication, provided that one uses the Chung-Hasan representation and
reduction methodology of \S\ref{sec:MRCPfieldrep}, to give a full modular multiplication algorithm.

Note that Karatsuba-Ofman multiplication~\cite{karatsuba} offers a similar trade-off for extension field arithmetic. 
Crucially however, as we show in \S\ref{subsec:costcomparison}, when we apply these formulae to GRPs the number of 
additions required is in fact reduced. One thus expects the CVMA to be significantly more efficient at contemporary ECC
bitlengths.
The original proof of~(\ref{eq:nogami_state}) given in~\cite{nogami} excludes some
intermediate steps and so for the sake of clarity we give a full proof in~\S\ref{subsec:derivation},
beginning with the following motivation.

\subsection{Alternative bases}\label{bases}

Observe that in the set of equations~(\ref{nogami_mul}), each of the
$2(m+1)$ coefficients $x_j,y_j$ is featured $m+1$ times, and so there is a
nice symmetry and balance to the formulae. However due
to the choice of basis, both $x_0$ and $y_0$ are implicitly assumed
to be zero. The output $\vec{z}$ naturally has this property also,
and indeed if one extends the multiplication algorithm to compute
$z_0$ we see that it equals $q_0 - q_0 = 0$.

At first sight, the expression $z_i = q_0 - q_i$ may seem a little unnatural.
It is easy to change the basis from $\{\omega^{m},\ldots,\omega\}$
to $\{\omega^{m-1},\ldots,\omega,1\}$:
for $\vec{x} = [x_{m-1},\ldots,x_{0}]$ and $\vec{y} = [y_{m-1},\ldots,y_{0}]$, we have:
\[
\vec{z} = \vec{x} \cdot \vec{y} = \sum_{i=0}^{m-1} z_i \omega^i,
\]
resulting in the expressions $z_i = q_m - q_i$,
with $q_i$ as given before. This change of basis relies on the relation
\begin{equation}\label{eliminate}
\omega^m \equiv - 1 - \omega - \cdots - \omega^{m-1} \bmod f(\omega).
\end{equation}
Note that in using this basis we have implicitly ensured that $x_m = y_m = 0$ in
~(\ref{nogami_mul}), rather than $x_0 = y_0 = 0$, 
and again the above formula is consistent since $z_m = q_m - q_m = 0$.
More generally if one excludes $\omega^k$ from the basis, then $x_k = y_k = 0$
and $z_i = q_k - q_i$.

One may infer from these observations that the most natural choice of basis
would seem to be $\{\omega^{m},\ldots,\omega,1\}$,
and that the expressions for $q_i$ arise from the arithmetic in the quotient 
ring $\mathcal{R}' = \F_p[\omega]/((\omega^{m+1}-1)\F_p[\omega])$, rather than $\F_{p^m} = \F_p[\omega]/(f(\omega)\F_p[\omega])$.
In this case multiplication becomes
\[
\vec{z} = \vec{x} \cdot \vec{y} = \sum_{i=0}^{m-1} z_i \omega^i = \sum_{i=0}^{m-1} (q_m - q_i) \omega^i =
\sum_{i=0}^{m} - q_i \omega^i,
\]
where for the last equality we have again used equation~(\ref{eliminate}).

\subsection{Derivation of coefficient formulae}\label{subsec:derivation}

We now derive the CVMA formulae of~(\ref{nogami_mul}).
Let $\vec{x}= [x_{m},\ldots,x_{0}] = \sum_{i=0}^{m} x_i \omega^i$, and
$\vec{y}= [y_{m},\ldots,y_{0}] = \sum_{i=0}^{m} y_i \omega^i$.
Then in the ring $\mathcal{R}'$, as in~(\ref{eq:convolution}) 
the product $\vec{x} \cdot \vec{y}$ is equal to $\sum_{i=0}^{m} z_i \omega^i$, where
\[
z_i = \sum_{j=0}^{m} x_{\langle j \rangle} y_{\langle i-j \rangle}.
\]
Of crucial importance is the following identity. For $0 \leq i \leq m$ we have:
\begin{equation}\label{identity}
2 \sum_{j=0}^{m} x_{\langle j \rangle} y_{\langle i-j \rangle}
- 2\sum_{j=0}^{m} x_{\langle j \rangle}y_{\langle j \rangle}
= - \sum_{j=0}^{m} (x_{\langle j \rangle} - x_{\langle i - j \rangle})
(y_{\langle j \rangle} - y_{\langle i - j\rangle}).
\end{equation}
To verify this identity observe that when one expands the terms in the right-hand side, the two
negative sums cancel with the second term on the left-hand side, since both are over a complete set of
residues modulo $m+1$. Similarly the two positive sums are equal and therefore cancel with the convolutions
in the first term on the left-hand side.
We now observe that there is some redundancy in the right-hand side of~(\ref{identity}),
in the following sense. First, observe that
\begin{equation}
\nonumber \sum_{j=0}^{m} x_{\langle \frac{i}{2} + j\rangle} y_{\langle \frac{i}{2} - j\rangle}
= \sum_{j=0}^{m} x_{\langle \frac{i}{2} + (j - \frac{i}{2}) \rangle} y_{\langle \frac{i}{2} - (j - \frac{i}{2}) \rangle}
= \sum_{j=0}^{m} x_{\langle j \rangle} y_{\langle i-j \rangle}.
\end{equation}
One can therefore rewrite the right-hand side of~(\ref{identity}) as:
\begin{equation}\label{trick}
- \sum_{j=0}^{m} (x_{\langle \frac{i}{2} + j\rangle} - x_{\langle \frac{i}{2} - j\rangle})
(y_{\langle \frac{i}{2} + j\rangle} - y_{\langle \frac{i}{2} - j\rangle}).
\end{equation}
Noting that the $j=0$ term of expression~(\ref{trick}) is zero, we rewrite it as:
\begin{equation}
\nonumber -\sum_{j=1}^{m/2} (x_{\langle \frac{i}{2} + j\rangle} - x_{\langle \frac{i}{2} - j\rangle})
(y_{\langle \frac{i}{2} + j\rangle} - y_{\langle \frac{i}{2} - j\rangle})
-\sum_{j=m/2+1}^{m} (x_{\langle \frac{i}{2} + j\rangle} - x_{\langle \frac{i}{2} - j\rangle})
(y_{\langle \frac{i}{2} + j\rangle} - y_{\langle \frac{i}{2} - j\rangle}),
\end{equation}
which in turn becomes
\begin{equation}
\nonumber -\sum_{j=1}^{m/2} (x_{\langle \frac{i}{2} + j\rangle} - x_{\langle \frac{i}{2} - j\rangle})
(y_{\langle \frac{i}{2} + j\rangle} - y_{\langle \frac{i}{2} - j\rangle})
-\sum_{j=1}^{m/2} (x_{\langle \frac{i}{2} - j\rangle} - x_{\langle \frac{i}{2} + j\rangle})
(y_{\langle \frac{i}{2} - j\rangle} - y_{\langle \frac{i}{2} + j\rangle}),
\end{equation}
and then upon negating the two terms in the second summation, we finally have
\begin{equation}
\nonumber - \sum_{j=0}^{m} (x_{\langle \frac{i}{2} + j\rangle} - x_{\langle \frac{i}{2} - j\rangle})
(y_{\langle \frac{i}{2} + j\rangle} - y_{\langle \frac{i}{2} - j\rangle}) =
2\sum_{j=1}^{m/2} (x_{\langle \frac{i}{2} + j\rangle} - x_{\langle \frac{i}{2} - j\rangle})
(y_{\langle \frac{i}{2} + j\rangle} - y_{\langle \frac{i}{2} - j\rangle}).
\end{equation}
 Hence~(\ref{identity}) becomes 
\begin{equation}\label{cyclicidentity}
\sum_{j=0}^{m} x_{\langle j \rangle} y_{\langle i-j \rangle} =
 \sum_{j=0}^{m} x_{\langle j \rangle}y_{\langle j \rangle}
 - \sum_{j=1}^{m/2} (x_{\langle \frac{i}{2} + j\rangle} - x_{\langle \frac{i}{2} - j\rangle})
(y_{\langle \frac{i}{2} + j\rangle} - y_{\langle \frac{i}{2} - j\rangle}).
\end{equation}
Equation~(\ref{cyclicidentity}) gives an expression for the
coefficients of the product $\vec{z}$ of elements $\vec{x}$ and
$\vec{y}$, in the ring $\mathcal{R}'$. Assuming these are computed using the more
efficient right-hand side, in order to restrict back to
$\F_p[\omega]/(f(\omega)\F_p[\omega])$, one can reduce the resulting polynomial
$\vec{z}$ by $f(\omega)$. Note however that one does not need to use a smaller
basis \`{a} la Nogami \ea in~\S\ref{subsec:mulnogami} or~\S\ref{bases} ,
but can reduce by $f(\omega)$ {\em implicitly}, 
without performing any computation. 
Indeed, letting $\langle \vec{x}, \vec{y} \rangle = \sum_{j=0}^{m} x_{\langle j \rangle}y_{\langle j\rangle}$, we have:
\begin{eqnarray}\label{cheap}
\nonumber \vec{z} &=& \sum_{i=0}^{m} z_i \omega^i = \sum_{i=0}^{m} 
(-q_i +\langle \vec{x}, \vec{y} \rangle) \omega^i = \sum_{i=0}^{m} -q_i\omega^i +
\langle \vec{x}, \vec{y} \rangle \sum_{i=0}^{m} \omega^i\\ 
&\equiv& \sum_{i=0}^{m} -q_i \omega^i \pmod{f(\omega)}.
\end{eqnarray}
Therefore the first term on the
right-hand side of~(\ref{cyclicidentity}) vanishes, so that one need not even compute it. 
Thus using the arithmetic in $\mathcal{R}'$ but implicitly working modulo $f(\omega)$
is more efficient than performing arithmetic in $\mathcal{R}'$ alone. This is somewhat fortuitous as
it means that while the multiply operation in~(\ref{cheap}) is not correct in $\mathcal{R}'$, nevertheless,
when one maps back to $\F_p[\omega]/(f(\omega)\F_p[\omega])$, it is correct.

\subsection{Application to GRPs}\label{subsec:application}

Since equation~(\ref{identity}) is an algebraic identity, it is easy to see that
exactly the same argument applies in the context of GRPs, and we can replace
the formulae~(\ref{eq:convolution}) with the CVMA formulae~(\ref{nogami_mul}). 
Since reduction in the ring $\mathcal{R} = \Z[T]/(T^{m+1}-1)\Z[T]$ has a particularly nice form for
GRPs, we choose to use the full basis for $\mathcal{R}$ and hence do not reduce {\em explicitly} modulo $\Phi_{m+1}(T)$ to obtain a
smaller basis. This also has the effect of eliminating the need to perform the addition of $q_0$ (or $q_m$, or
whichever term one wants to eliminate when one reduces modulo $\Phi_{m+1}(T)$), 
simplifying the multiplication algorithm further.
Absorbing the minus sign into the $q_i$, Algorithm~\ref{alg:mrcpmul}
details how to multiply residue representatives. 

\begin{remark} 
Observe that each component of $\vec{z}$ may be computed entirely
independently of the others. Hence using $m+1$ processors rather than $1$, it
would be possible to speed up the execution time of
Algorithm~\ref{alg:mrcpmul} by a factor of $m+1$, making it particularly
suitable for hardware implementation. In~\S\ref{sec:reduction} we consider the
parallelisation of our reduction algorithms as well.
\end{remark}

\begin{algorithm}{GRP MULTIPLICATION}
{$\vec{x} = [x_{m},\ldots,x_{0}], \vec{y} = [y_{m},\ldots,y_{0}]
\in \Z^{m+1}$} 
{$\vec{z} = [z_{m},\ldots,z_{0}] \in \Z^{m+1}$ \\
\hspace{15mm} where $\vec{z} \cong_{\Phi_{m+1}(t)} \vec{x} \cdot \vec{y}$}\label{alg:mrcpmul}
\nline For $i=m$ to $0$ do: \\
\nnline
\hspace{-5mm} $z_i \leftarrow \sum_{j=1}^{m/2} (x_{\langle \frac{i}{2} - j\rangle} - x_{\langle \frac{i}{2} + j\rangle})\cdot
(y_{\langle \frac{i}{2} + j\rangle} - y_{\langle \frac{i}{2} - j\rangle})$\\
\nline Return $\vec{z}$
\end{algorithm}

\subsection{Cost comparison}\label{subsec:costcomparison}

We here use a simple cost model to provide a measure of the potential performance 
improvement achieved by using Algorithm~\ref{alg:mrcpmul}, rather 
than schoolbook multiplication of residues. We assume the inputs to the multiplication algorithm have coefficients 
bounded by $b^q$, \ie they each consist of $q$ words.
Let $M(q,q)$ be the cost of a $q$-word by $q$-word
schoolbook multiplication, and let $A(q,q)$ be the cost of an ition of two
$q$-word values. We assume that $A(2q,2q) = 2A(q,q)$ and that there is
no overflow beyond $2q$ words in the resulting vector components, which one
can ensure by selecting appropriate GRPs, see~\S\ref{sec:stability}.
The cost of the multiplication using each method is as follows.

\subsubsection{GRP schoolbook multiplication}

Working modulo $T^m + \cdots + T + 1$ and using a basis consisting of $m$
terms only, the number of coefficient multiplications is $m^2$, while the
number of double-length additions is also $m^2$. Hence the total cost is simply
\[
m^2 \cdot M(q,q) + 2m^2 \cdot A(q,q).
\]
Note that computing the convolution~(\ref{eq:convolution}) costs 
\[(m+1)^2 \cdot M(q,q) + 2m(m+1) \cdot A(q,q),
\]
which is costlier since it requires embedding into $\mathcal{R}$, which introduces some redundancy.

\subsubsection{CVMA formulae}

For each $z_i$ computing each term in the sum costs $M(q,q) + 2A(q,q)$,
and so computing all these terms costs $\frac{m}{2} \cdot (M(q,q) + 2A(q,q))$.
The cost of adding these is $(\frac{m}{2} - 1) A(2q,2q) = (m - 2)\cdot A(q,q)$.
For all the $m+1$ terms $z_i$ the total cost is therefore
\[
\frac{m(m+1)}{2} \cdot M(q,q) + 2(m^2-1)\cdot A(q,q).
\]
Therefore by using the CVMA formulae, we reduce not only the number of multiplications,
but also the number of additions (by $2$), contrary to the case of field
extensions, for which the CVMA formulae increases the number of additions by
nearly $50\%$. We have thus found an analogue of the asymptotic cyclic versus linear
convolution speed-up for multiplication modulo Mersenne numbers 
(see Eq. (6.1) of~\cite{crandallfagin}, for example) at small bitlengths for which schoolbook multiplication
is optimal, for GRPs.

%------------------------------------------------------------------------------------------------
%------------------------------------------------------------------------------------------------

\section{GRP Reduction}
\label{sec:reduction}

In this section we detail reduction algorithms for two types of GRPs.
The first, Algorithm~\ref{alg:red1}, assumes only that $t$ is even, which
provides the minimum possible restriction on the form of the resulting GRPs
for any given bitlength. All such GRPs can therefore be implemented 
with code parametrised by the single variable $t$, which may be beneficial for
some applications. Supposing that $R = b^q > t$, then as with Montgomery reduction, it is more efficient to reduce 
components not by $R$ as in~(\ref{basic_mont1}) and~(\ref{basic_mont2}), but
by $b$ sequentially $q$ times. In Algorithm~\ref{alg:red1} each reduction 
therefore reduces the input's components by approximately $\log_2{b}$ bits.

The second reduction method as detailed in Algorithm~\ref{alg:red2} is a
specialisation of Algorithm~\ref{alg:red1}. It assumes that $t \equiv 0 \bmod 2^l$ for some $ l > 1$,
and each application of the reduction function reduces the input's components by approximately
$l$ bits. Algorithm~\ref{alg:red2} is potentially far more efficient than Algorithm~\ref{alg:red1},
depending on the form of $t$. Ideally one should choose a $t$ for which $l > (\log_2{t})/2$ so that
two applications of the reduction function are sufficient in order to produce
components of the desired size, which is minimal. In general for other values
of $l$ a larger number of reductions may be needed, which we consider in \S\ref{sec:stability}.
In constrast to Algorithm~\ref{alg:red1}, which is designed for
generality, Algorithm~\ref{alg:red2} is geared towards high-speed reduction. The trade-off arising
here is that there will naturally be far fewer GRPs of this restricted form.
We also present a modification of Algorithm~\ref{alg:red2}, which is slightly
more efficient in practice, in Algorithm~\ref{alg:red3}.

\subsection{GRP reduction: $t$ even}

Following \S\ref{subsec:montlattice}, in equation~(\ref{basic_mont1}) we need
the matrix $-L^{-1}$:
\begin{equation}
-L^{-1} = \frac{1}{t^{m+1}-1}\left[ \begin{array}{cccccc}
1  & t^{m} & \cdots & t^3 & t^2 & t \\
t & 1 & \cdots & t^4 & t^3 & t^2   \\
t^2 & t & \cdots & t^5 & t^4 & t^3   \\
\vdots & \vdots & \ddots & \vdots &\vdots & \vdots \\
t^{m-1} & t^{m-2} & \cdots & t & 1 & t^{m} \\
t^{m} & t^{m-1} & \cdots & t^2 & t & 1
\end{array} \right].
\end{equation}
The form of $L$ and $-L^{-1}$ allows one to compute $\vec{u} = -L^{-1} \cdot \vec{z} \bmod b$ and 
$L \cdot \vec{u}$, computed in equation~(\ref{basic_mont2}), very efficiently. Since $t$ is even, the following
vector may be computed. Let $t[0]$ be the least significant digit of $t$, written in base $b$,
and let
\[
\vec{V} \stackrel{\rm def}{=} \frac{1}{t[0]^{m+1}-1}[t[0]^{m}, t[0]^{m-1},\ldots,t[0],1] \bmod b.
\] 
Algorithm~\ref{alg:red1} details how to reduce a given an input vector $\vec{z}$ by $b$, modulo $t^{m+1}-1$, 
given the precomputed vector $\vec{V}$. Observe that Algorithm~\ref{alg:red1} greatly simplifies 
the reduction algorithm originally given in~\cite{CH3}.
This is possible since for $t^{m+1}-1$ one can interleave the computation 
of the vectors $\vec{u}$ and $\vec{w}$ defined in~(\ref{basic_mont1})
and~(\ref{basic_mont2}) respectively. This has two benefits. First, as one
computes each component of $\vec{w}$ sequentially, one need only store a single component of $\vec{u}$, rather than $m+1$.
Second, since when one computes $L \cdot \vec{u}$ one needs to compute $t\cdot u_{\langle i+1 \rangle}$
for $i=m,\ldots,0$ (in \url{line} \url{3}), one obtains $t[0] \cdot u_i$ (the
first term on right-hand side of \url{line} \url{4}) 
for free by computing the full product $t \cdot u_{\langle i+1 \rangle}$ first. 
One therefore avoids recomputing the least significant digit of $t \cdot u_{\langle i+1 \rangle}$ in each loop iteration.
In fact one can do this for
\begin{algorithm}{$\mathrm{red1}_b(\vec{z})$}
{$\vec{z} = [z_{m},\ldots,z_{0}] \in \Z^{m+1}$}
{$\mathrm{red}_{b}(\vec{z})$ where $\mathrm{red}_{b}(\vec{z}) \cong_{t^{m+1}-1} \vec{z} \cdot b^{-1}$}
\label{alg:red1}
\nline Set $u_0 \leftarrow (\sum_{i=0}^{m} V_i \cdot z_{i}[0]) \bmod b$\\
\nline For $i=m$ to $0$ do: \\
\nnline $v_i \leftarrow t \cdot u_{\langle i+1 \rangle}$\\
\nnline $u_i \leftarrow (v_i[0] - z_{i}[0]) \bmod b$\\
\nnline $w_i \leftarrow (z_i + u_i - v_i)/b$\\
\nline Return $\vec{w}$
\end{algorithm}
\noindent any polynomial $t^{m+1} - c$, with exactly the same algorithm, the only difference being
in the definition of $\vec{V}$, where $t^{m+1}-c$ becomes the denominator. 
For polynomials with other non-zero coefficients, this does not seem possible,
and so Algorithm~\ref{alg:red1} seems likely to be the most efficient Chung-Hasan
reduction possible with this minimal restriction on the form of $t$.

It is straightforward to verify that Algorithm~\ref{alg:red1} correctly produces an output vector 
in the correct congruency class, via a sequence of simple transformations
of~\cite[Algorithm 3]{CH3}. However we do not do so here, since we are mainly interested in the more
efficient Algorithms~\ref{alg:red2} and~\ref{alg:red3}.

\begin{remark}
Note that in the final loop iteration, $u_0$ from \url{line} \url{1} is recomputed,
which is therefore unnecessary. However, we chose to write the algorithm in this form to emphasise its cyclic
structure. Indeed, there is no need to compute $u_0$ first;
if one cyclically rotates $\vec{V}$ by $j$ places to the left, then the
vector $\vec{w}$ to be added to $\vec{z}$ in~(\ref{basic_mont2}) is rotated $j$ places to the left also.
One can therefore compute each coefficient of
$\mathrm{red1}_{b}(\vec{z})$ independently of the others using a 
rotated definition for $\vec{V}$ (or equivalently by rotating the input $\vec{z}$ ).
This demonstrates that a parallelised version of the reduction algorithm with
$m+1$ processors is feasible. However, as each processor requires the least
significant word of each component of $\vec{z}$, this necessitates a synchronised
broadcast before each invocation of the reduction function. In
this scenario the reduction time would be proportional to the number of such
broadcasts and reductions required, independently of $m+1$. 
\end{remark}

\subsection{GRP reduction: $t \equiv 0 \bmod{2^l}$}
\label{sec:redfast}

In the ideal case that $t=2^l$, we see that such a GRP would be a GMN. 
In this case, one can use the reduction method detailed in \S\ref{subsec:latticereduction}
without resorting to using its Montgomery version at all. 
Multiplication would also be faster thanks to Nogami's formulae.
Unfortunately, such GRPs seem to be very rare. It is easy to show that if $t =
2^l$ with $l > 1$ and $\Phi_{m+1}(t)$ is prime, then $l = m+1$. Testing the first few cases, we
find prime GRPs for $l = 2,3,7,59$ but no others for prime $l < 400$. 
Note that these primes contradict Dubner's assertion that no such GRPs exist~\cite[\S2]{dubner3}.
Since for $l = 59$ the corresponding GRP has $3422$ bits, this is
already out of our target range for ECC, so we need not worry about such GRPs.

Hoping not to cause confusion, in this subsection we now let $b=2^l$ 
where $l$ is not necessarily and usually not the
word size of the target architecture. We denote the cofactor of $b$ in $t$
by $c$ (which by the above discussion we assume is $>1$), so that $t = b \cdot
c$. Algorithm~\ref{alg:red2} details how to reduce a given an input vector $\vec{z}$ by $b$, modulo $t^{m+1}-1$.
\begin{algorithm}{$\mathrm{red2}_b(\vec{z})$}
{$\vec{z} = [z_{m},\ldots,z_{0}] \in \Z^{m+1}$}
{$\mathrm{red}_{b}(\vec{z})$ where $\mathrm{red}_{b}(\vec{z})\cong_{t^{m+1}-1} \vec{z} \cdot b^{-1}$} \label{alg:red2}
\nline For $i=m$ to $0$ do: \\
\nnline \hspace{-5mm} $w_i \leftarrow (z_i + (-z_i \bmod{b}))/b - c\cdot(-z_{\langle i+1 \rangle} \bmod{b})$\\
\nline Return $\vec{w}$
\end{algorithm}
A simple proof of correctness of Algorithm~\ref{alg:red2} comes from the specialisation of Algorithm~\ref{alg:red1}.
Since $t \equiv 0 \bmod{b}$, writing $t$ in base $b$, the vector $\vec{V}$ becomes
\[
\vec{V} \stackrel{\rm def}{=} [0,\ldots,0 ,-1] \bmod b.
\] 
Hence for \url{line} \url{1} of Algorithm~\ref{alg:red1} we have 
\[
u_0 \leftarrow -z_0[0] \bmod{b}.
\]
Since in \url{line} \url{4} of Algorithm~\ref{alg:red1}, we have $v_i \equiv 0 \bmod{b}$, we deduce that $u_i = - z_i \bmod{b}$, 
and hence we can eliminate $u_i$ altogether. Each loop iteration then simplifies to
\begin{eqnarray}
\nonumber &v_i& \leftarrow t \cdot ( - z_{\langle i+1 \rangle} \bmod{b})\\
\label{linez} &w_i& \leftarrow (z_i + (-z_i \bmod{b}) - v_i)/b
\end{eqnarray}
Upon expanding~(\ref{linez}), we obtain
\begin{eqnarray}
\nonumber w_i &\leftarrow& (z_i + (-z_i \bmod{b}))/b - t \cdot ( - z_{\langle i+1 \rangle} \bmod{b})/b \\
\nonumber     &= & (z_i + (-z_i \bmod{b}))/b - c \cdot ( - z_{\langle i+1 \rangle} \bmod{b}),
\end{eqnarray}
as required. However since we did not provide a proof of correctness of Algorithm~\ref{alg:red1}, we also give a direct proof as follows.
Observe that modulo $t^{m+1} - 1$, we have 
\begin{eqnarray}
\nonumber \psi^{-1}(\vec{w}) &\equiv& \sum_{i=0}^{m} w_i t^i\\
\nonumber &\equiv& \sum_{i=0}^{m} [(z_i + (-z_i \bmod{b}))/b - c \cdot(-z_{\langle i+1 \rangle} \bmod{b})]t^i\\
\nonumber &\equiv& \sum_{i=0}^{m} (z_i/b)t^i  + \sum_{i=0}^{m} (-z_i \bmod{b}))/b)t^i - \sum_{i=0}^{m} ((-z_{\langle i+1 \rangle} \bmod{b})/b)t^{i+1}\\
\nonumber &\equiv& \sum_{i=0}^{m}z_it^i/b \pmod{t^{m+1}-1}
\end{eqnarray}
as required. In terms of operations that may be performed very efficiently, we alter Algorithm~\ref{alg:red2} slightly to give 
Algorithm~\ref{alg:red3}, which has virtually the same proof of correctness as the one just given.
\begin{algorithm}{$\mathrm{red3}_b(\vec{z})$}
{$\vec{z} = [z_{m},\ldots,z_{0}] \in \Z^{m+1}$}
{$\mathrm{red}_{b}(\vec{z})$ where $\mathrm{red}_{b}(\vec{z})\cong_{t^{m+1}-1} \vec{z} \cdot b^{-1}$} \label{alg:red3}
\nline For $i=m$ to $0$ do: \\
\nnline \hspace{-5mm} $w_i \leftarrow z_i/b + c\cdot(z_{\langle i+1 \rangle} \bmod{b})$\\
\nline Return $\vec{w}$
\end{algorithm}
Note that the first term in \url{line} \url{2} of Algorithm~\ref{alg:red2} has been replaced by a division by $b$, which 
can be effected as a simple shift, while now the second term needs the positive residue modulo $b$,
which can be extracted more efficiently. Hence Algorithm~\ref{alg:red3} is the one we use.
By our previous discussion, $c$ necessarily has Hamming weight at least two
for GRPs in our desired range. By using $c$ that have very low Hamming weight,
one can effect the multiplication by $c$ by shifts and adds, rather than a
multiply (or \url{imulq}) instruction. Hence for such GRPs, assuming only two
invocations of Algorithm~\ref{alg:red3} are needed, reduction will be
extremely efficient.

\begin{remark}
Regarding parallelisation, observe that for $m+1$ processors, only
the least significant word of $z_{\langle i+1 \rangle}$ is passed to processor
$i$, thus reducing the broadcast requirement in comparison with Algorithm~\ref{alg:red1}. 
\end{remark}

%-----------------------------------------------------------------------------------------------------
%-----------------------------------------------------------------------------------------------------

\section{GRP Residue Representation}
\label{sec:representation}

So far in our treatment of both multiplication and reduction, for the sake of generality 
we have assumed arbitrary precision when representing GRP residues in $\Z^{m+1}$. 
In this section we specialise to fixed precision and develop a residue representation that
ensures that our chosen algorithms are efficient.  
Our decisions are informed purely by our chosen multiplication and reduction algorithms 
--- Algorithms~\ref{alg:mrcpmul} and~\ref{alg:red3} ---
which we believe offer the best performance for GRPs for the relatively small bitlengths
which are relevant to ECC.
In other scenarios or if considering asymptotic performance, one would need to redesign the residue representation and 
multiplication algorithm accordingly. 

For $x \in \{0,\ldots,t^{m+1} - 1\}$ we write $\vec{x} = [x_m,\ldots,x_0]$
for its base-$t$ expansion, \ie  $x = \sum_{i=0}^{m} x_i t^i$.
The base-$t$ representation has positive coefficients, 
however Algorithm~\ref{alg:mrcpmul} makes use of negative coefficients, so we prefer to incorporate these.
We therefore replace the mod function in the conversion with mods, 
the least absolute residue function, to obtain a residue in the interval
$[-t/2,t/2-1]$:
\[
\mathrm{mods}(x) =
\left\{
    \begin{array}{ll}
      x \bmod{t} &           \hspace{3mm} \mathrm{if} \ (x \bmod{t}) < t/2, \\
       x \bmod{t} - t      & \hspace{3mm} \mathrm{otherwise}.
    \end{array}
\right.
\]
Using this function, Algorithm~\ref{alg:map2poly} converts residues modulo $t^{m+1} - 1$ into
the required form~\cite[Algorithm 1]{CH3}.
\begin{algorithm}{BASE-$t$ CONVERSION $\psi$}
{An integer $0 \leq x < t^{m+1} - 1$}
{$\vec{x} = [x_m.\ldots,x_0]$ such that $|x_i| \leq t/2$\\ 
\hspace{14.5mm} and $\sum_{i=0}^{m} x_i t^i \equiv x \pmod{t^{m+1}-1}$}\label{alg:map2poly}\\
\nline For $i$ from $0$ to $m$ do: \\
\nnline $x_i \leftarrow x \ \mathrm{mods} \ t$\\
\nnline $x \leftarrow (x - x_i)/t$\\
\nline $x_0 \leftarrow x_0 + x$\\
\nline Return $\vec{x} = [x_m,\ldots,x_0]$
\end{algorithm}

The reason for \url{line} \url{4} in Algorithm~\ref{alg:map2poly} is to reduce modulo $t^{m+1} - 1$ the 
coefficient of $t^{m+1}$ possibly arising in the expansion. Note that in this addition, $x \in \{0,1\}$, 
and hence $|x_i| \leq t/2$ for each $0 \leq i \leq m$. By construction, we in fact 
have $-t/2 \leq x_i < t/2$ for $1 < i < m$ while only $x_0$ can attain the upper bound of $t/2$.
There are therefore $t^{m}(t + 1)$ representatives in this format, thus
introducing a very small additional redundancy.
Letting $k = \lceil \log_{2} t \rceil$, if we assume $t \leq 2^k - 2$, so that $[-t/2,t/2] \subset [-2^k/2,2^k/2 - 1]$,
then the coefficients as computed above can be represented in two's complement in $k$ bits.
In terms of efficiency, Algorithm~\ref{alg:map2poly} contains divisions by
$t$, which requires not only time, but also space, which on some
platforms may be at a premium. Writing $t = 2^l \cdot c$ as in
\S\ref{sec:redfast}, then if the cofactor $c = 2^{k-l} - c'$ with $c'$ very
small, then division by $t$ consists of a shift right by $l$ bits and a
division by $c$, which can be performed efficiently using Algorithm 1 of~\cite{CH1}.

Following this conversion, it might seem desirable to define vectors whose components are in $[-2^k/2,2^k/2 - 1]$
to be reduced, or canonical residue representatives. However,
for efficiency purposes it is preferable to have a reduction function which,
when performed sufficiently many times, outputs an element for which
one does not have to perform any modular additions or subtractions to make
reduced, as this eliminates data-dependent branching.
A control-flow invariant reduction function is also essential to defend
against side-channel attacks, see \S\ref{sec:otherops}.
To obtain such a function, observe that the second term in \url{line} \url{2} of
Algorithm~\ref{alg:red3}, namely $c\cdot (z_{\langle i+1 \rangle} \bmod{b})$, 
is positive, and in the worst case is $k$ bits long. The first term, $z_i/b$,
is clearly $l = \log_2 b$ bits shorter than $z_i$. Since one adds these the resulting value may be $k+1$ bits, or larger,
depending on the initial length of the inputs' components. Furthermore, since we wish to allow negative 
components, in two's complement the output requires a further bit, giving a minimal requirement of $k+2$ bits.
We therefore choose not to use minimally reduced elements as coset representatives in 
$\Z^{m+1}/\sim$, as output by Algorithm~\ref{alg:map2poly}, but slightly larger elements, which we now define.
\begin{definition}
We define the following set of elements of $\Z^{m+1}$ to be {\em reduced}:
\begin{equation}\label{reduced}
\I^{m+1} = \{ [x_m,\ldots,x_0] \in \Z^{m+1} \mid - 2^{k+1} \leq x_i < 2^{k+1}\}.
\end{equation}
\end{definition}
Note that the redundancy inherent in this representation depends on how close $t$ is to $2^{k+2}$.
For a modular multiplication, we assume that the inputs are reduced. We must therefore ensure that
the output is reduced also. This naturally leads one to consider I/O stability, as we do in \S\ref{sec:stability}.

Once we have a reduced representative $\vec{x} = \psi(x)$ we also need to convert to the Montgomery domain. 
While one can do this in $\Z/(t^{m+1}-1)\Z$ before applying $\psi$, it is more convenient to do so in $\Z^{m+1}/\sim$. 
Assuming $q$ reductions by $b$ are sufficient to ensure I/O modular multiplication stability, 
we precompute $\psi(b^{2q} \bmod{\Phi_{m+1}(t)})$ and then using Algorithms~\ref{alg:mrcpmul} and~\ref{alg:red3} compute 
\[
\vec{x} \cdot \psi(b^{2q} \bmod{\Phi_{m+1}(t)})/b^q \cong_{\Phi_{m+1}(t)} \psi(x \cdot b^q).
\]
Similarly, to get back from the Montgomery domain, again using Algorithms~\ref{alg:mrcpmul} and~\ref{alg:red3}, we compute
\[
\psi(x \cdot b^q) \cdot \psi(1) /b^q \cong_{\Phi_{m+1}(t)} \psi(x).
\]
With regard to mapping back from $\vec{x} = [x_m,\ldots,x_0] \in \I^{m+1}$ to canonical residues in $\Z/\Phi_{m+1}(t)\Z$, one has 
\[
\sum_{i=0}^{m} x_i t^i \equiv \sum_{i=0}^{m-1} (x_i - x_m)t^i \pmod{\Phi_{m+1}(t)},
\]
which can be computed efficiently by first using Horner's rule and then mapped to $\{0,\ldots,\Phi_{m+1}(t)-1\}$ 
by repeated additions or subtractions. In terms of operations required for ECC, we assume that the conversions 
are one-time computations only, with all other operations taking place in the (Montgomery) Chung-Hasan representation.

%-----------------------------------------------------------------------------------------------------
%-----------------------------------------------------------------------------------------------------

\section{Modular Multiplication Stability}
\label{sec:stability}

In this section we analyse Algorithms~\ref{alg:mrcpmul} and~\ref{alg:red3}
with a view to ensuring I/O stability for modular multiplication. 
We assume the following: $b=2^l$, $t=c \cdot b$ where $c < 2^{k-l}$ (and hence $t < 2^k - 2$),
and that reduced elements have the form~(\ref{reduced}).
Input elements therefore have components in $\I = [-2^{k+1},2^{k+1} - 1]$, and these are representable in $k+2$ bits in
two's complement. For simplicity and in order for our analysis to be as general as possible, we use the term single precision
to mean a word base large enough to contain $t$ --- even if this in fact requires multiprecision on a given architecture ---
and double precision to mean twice this size. We assume that for this single precision word size $w$,
the components of $\vec{z}$ output by Algorithm~\ref{alg:mrcpmul} are double precision.
In practice one prefers to specialise to actual single precision $t$ on a given architecture, 
since this obviates the need for multiprecision arithmetic; utilising the native double 
precision multipliers that most CPUs possess is more efficient, and reduction is also faster 
for smaller $t$ since fewer iterations need be performed. 
We note that in constrained environments however, multiprecision may however be unavoidable.

During the multiplication, terms of the form $x_i - x_j$ are computed, which are bounded by
\[
-2^{k+2} + 1 \leq x_i - x_j \leq 2^{k+2} - 1,
\]
and which therefore fit into $k+3$ bits in two's complement. The product of two such elements is performed,
giving a result
\[
-2^{2k+4} + 2^{k+3} - 1 \leq (x_i - x_j)\cdot(y_j-y_i) \leq 2^{2k+4} - 2^{k+3} + 1,
\]
which fits into $2k+5$ bits in two's complement. One then adds $m/2$ of these terms, giving a possible expansion
of up to $\lceil \log_2{m/2} \rceil$ bits, which must be double precision. 
We therefore have a constraint on the size of $t$ (in addition to the constraint $t < 2^k - 2$) in terms of $m$:
\begin{equation}\label{constraint}
\lceil \log_2{(m/2)} \rceil + 2k + 5 \leq 2w
\end{equation}
This inequality determines a constraint on the size of $t$, given $m$ and $w$. 
Assuming~(\ref{constraint}) is satisfied, one then needs to find the minimum value of $b = 2^l$ such that the result
of the multiplication step, when reduced by $b$ a specified number of times, say $q$, outputs a reduced element. 
This needs to be done for each $(m,k)$ found in the procedure above. 
Any power of $2$ larger than this minimum will obviously be satisfactory also, however minimising $b$ maximises the set of 
prime-producing cofactors $c$, which as stated in \S\ref{sec:reduction} may be
useful in some scenarios.

In~\S\ref{sec:representation}, we showed that one application of Algorithm~\ref{alg:red3}
shortened an input's components by $l-1$ bits, unless the components were already shorter than  
$(k+2) + (l-1)$ bits. Therefore stipulating that $q$ reductions suffice to produce a reduced output,
we obtain a bound on $l$ in the following manner. Let 
\[
h = \lceil \log_2{(m/2)} \rceil + 2k + 5
\]
Then after one reduction, the maximum length of a component is $h - l + 1$. Similarly
after $q$ reductions, the maximum length is $\max\{h - q (l-1), k+2    \}$, 
and we need this to be at most $k+2$. Hence our desired condition is
\[
h - q(l-1) \leq k+2
\]
Solving for $l$, we have
\begin{equation}\label{eq:lbound}
l \geq 1 + \frac{\lceil \log_2{(m/2)} \rceil + k + 3}{q}
\end{equation}

Using these inequalities it is an easy matter to generate triples $(m+1,k,l)$ which ensure multiplication
stability for any $w$ and $q$.
For example, for $w=64$, Tables~\ref{stableparameters1} and~\ref{stableparameters2}
give sets of stable parameters for $q=2$ and $q=3$ respectively.
\begin{table}\caption{Stable parameters: $w=64$, $q=2$}\label{stableparameters1}
\begin{center}
\begin{tabular}{|c|c|c|c|c|}
\hline
$m+1$ & $k$ & $l$ & $c <$ & $\lceil \log_2{p} \rceil $\\
\hline
$3$   & $61$   &  $33$   &   $2^{28}$ &    $122$\\ 
$5$   & $61$   &  $34$   &   $2^{27}$ &    $244$\\ 
$7$   & $60$   &  $34$   &   $2^{26}$ &    $360$\\ 
$11$  & $60$   &  $34$   &   $2^{26}$ &    $600$\\ 
$13$  & $60$   &  $34$   &   $2^{26}$ &    $720$\\ 
$17$  & $60$   &  $34$   &   $2^{26}$ &    $960$\\ 
\hline
\end{tabular}
\end{center}
\end{table}

\begin{table}\caption{Stable parameters: $w=64$, $q=3$}\label{stableparameters2}
\begin{center}
\begin{tabular}{|c|c|c|c|c|}
\hline
$m+1$ & $k$ & $l$ & $c <$ & $\lceil \log_2{p} \rceil $\\
\hline
$3$   & $61$   &  $23$   &   $2^{38}$ &    $122$\\ 
$5$   & $61$   &  $23$   &   $2^{38}$ &    $244$\\ 
$7$   & $60$   &  $23$   &   $2^{37}$ &    $360$\\ 
$11$  & $60$   &  $23$   &   $2^{37}$ &    $600$\\ 
$13$  & $60$   &  $23$   &   $2^{37}$ &    $720$\\ 
$17$  & $60$   &  $23$   &   $2^{37}$ &    $960$\\ 
\hline
\end{tabular}
\end{center}
\end{table}

The final column gives the maximum bitlength of a GRP that can be represented with those parameters,
though of course by using smaller $c$ one can opt for smaller primes, and the
corresponding minimum value of $l$ reduces according to~(\ref{eq:lbound}).
To generate suitable GRPs, a simple linear search over the values of $c$ of the desired size is sufficient, 
checking whether or not $\Phi_{m+1}(2^l \cdot c)$ is prime, see~\S\ref{sec:paramgen}.

%---------------------------------------------------------------------------------------------------
%---------------------------------------------------------------------------------------------------

\section{Full GRP Modular Multiplication}
\label{sec:fullmul}

For completeness we now piece together the parts treated thus far into a full modular
multiplication algorithm, where in Algorithm~\ref{modmul} we assume $q$
reductions by $b$ are required for I/O stability and in \url{line} \url{4} either 
Algorithm~\ref{alg:red1} or Algorithm~\ref{alg:red3} is used according to
the form of $b$.

\begin{algorithm}{GRP MODMUL}
{$\vec{x} = [x_{m},\ldots,x_{0}], \vec{y} = [y_{m},\ldots,y_{0}]
\in \I^{m+1}$} {$\vec{z} = [z_{m},\ldots,z_{0}] \in \I^{m+1}$ where 
$\vec{z} \cong_{\Phi_{m+1}(t)} \vec{x} \cdot \vec{y} \cdot b^{-q}$}\label{modmul}
\nline For $i=m$ to $0$ do: \\
\nnline
\hspace{-5mm} $z_i \leftarrow \sum_{j=1}^{m/2} (x_{\langle \frac{i}{2} - j\rangle} - x_{\langle \frac{i}{2} + j\rangle})\cdot
(y_{\langle \frac{i}{2} + j\rangle} - y_{\langle \frac{i}{2} - j\rangle})$\\
\nline For $k$ from $0$ to $q-1$ do: \\
\nnline \hspace{-5mm} $\vec{z} \leftarrow \mathrm{red}_b(\vec{z})$ \\
\nline Return $\vec{z}$
\end{algorithm}

Should $t$ be multiprecision on a particular architecture, then as with Montgomery
arithmetic it may be more efficient to use an interleaved multiplication and reduction
algorithm, as we detail in Algorithm~\ref{interleaved}. Here one needs $b$ to be the word base
of the underlying architecture and so in \url{line} \url{6}, if $t \equiv 0 \pmod{b}$ we use 
Algorithm~\ref{alg:red3}, otherwise we use Algorithm~\ref{alg:red1}. For 
$\vec{x} = [x_m,\ldots,x_0]$ we write $x_i = x_i[0] + x_i[1]b + \cdots + x_i[q-1]b^{q-1}$.  

\begin{algorithm}{GRP MODMUL (interleaved)}
{$\vec{x} = [x_{m},\ldots,x_{0}], \vec{y} = [y_{m},\ldots,y_{0}] \in \I^{m+1}$}
{$\vec{z} = [z_{m},\ldots,z_{0}] \in \I^{m+1}$ where
$\vec{z} \cong_{\Phi_{m+1}(t)} \vec{x} \cdot \vec{y} \cdot b^{-q}$}
\label{interleaved}
\nline $\vec{z} \leftarrow [0,\ldots,0]$ \\
\nline For $k=0$ to $q-1$ do: \\
\nnline For $i=m$ to $0$ do: \\
\nnnline $w_{i} \leftarrow \sum_{j=1}^{m/2} (x_{\langle \frac{i}{2} -j \rangle}[k] - x_{\langle \frac{i}{2} + j \rangle}[k]) \cdot 
(y_{\langle \frac{i}{2} + j \rangle} - y_{\langle \frac{i}{2} - j \rangle})$ \\
\nnline $\vec{z} \leftarrow \vec{z} + \vec{w}$ \\
\nnline $\vec{z} \leftarrow \text{red}_{b}(\vec{z})$\\
\nline Return $\vec{z}$
\end{algorithm}

To verify the correctness of Algorithm~\ref{interleaved}, observe that
for each of the $m+1$ components of $\vec{z}$, after the last iteration of the outer loop we have: 
\begin{eqnarray}
\nonumber z_i &=& \sum_{j=1}^{m/2} \Big(\sum_{k=0}^{q-1} (x_{\langle
  \frac{i}{2} -j \rangle}[k] - x_{\langle \frac{i}{2} + j \rangle}[k])/b^{q-k}\Big) \cdot (y_{\langle \frac{i}{2} + j \rangle} - y_{\langle \frac{i}{2} - j \rangle})\\
\nonumber  &=& \sum_{j=1}^{m/2} ((x_{\langle \frac{i}{2} -j \rangle} - x_{\langle \frac{i}{2} + j \rangle})/b^{q}) \cdot
(y_{\langle \frac{i}{2} + j \rangle} - y_{\langle \frac{i}{2} - j \rangle}).
\end{eqnarray}
Hence when taken modulo $\Phi_{m+1}(t)$, we see that $\vec{z}$ is congruent to:
\begin{eqnarray}
\nonumber \sum_{i=0}^{m} z_i \cdot t^i
&\cong_{\Phi_{m+1}(t)}& \sum_{i=0}^{m} \Big(\sum_{j=1}^{m/2} (x_{\langle \frac{i}{2} -j \rangle} - x_{\langle \frac{i}{2} + j \rangle})/b^q) 
\cdot (y_{\langle \frac{i}{2} + j \rangle} - y_{\langle \frac{i}{2} - j \rangle})\Big) \cdot t^i\\
\nonumber &\cong_{\Phi_{m+1}(t)}& \sum_{i=0}^{m} \Big( \sum_{j=1}^{m/2} (x_{\langle \frac{i}{2} -j \rangle} - x_{\langle \frac{i}{2} + j \rangle}) 
\cdot (y_{\langle \frac{i}{2} + j \rangle} - y_{\langle \frac{i}{2} - j \rangle}) \cdot t^i \Big) /b^q\\ 
\nonumber &\cong_{\Phi_{m+1}(t)}& \vec{x} \cdot \vec{y} \cdot b^{-q}, 
\end{eqnarray}
as required. As with ordinary Montgomery arithmetic, there are many possible
ways to perform the interleaving, see~\cite{koc} for example.

%--------------------------------------------------------------------------------------------
%--------------------------------------------------------------------------------------------

\section{Other arithmetic and side-channel secure ECC}\label{sec:otherops}

In addition to modular multiplication, one also needs to perform other arithmetic operations when
implementing ECC point multiplication. In this section we detail how to perform these using our 
representation and briefly explain how it enables point multiplication to be made immune to various side-channel attacks.

\subsection{Other arithmetic operations}

\subsubsection{Addition/subtraction}
\label{sec:addition}

To perform an addition or subtraction of two reduced elements $\vec{x},\vec{y}$, we compute the following:
\[
\vec{x} \pm \vec{y} = [x_m \pm y_m,\ldots,x_0 \pm y_0]. 
\]
Note that the bounds on each of these components is $[-2^{k+2},2^{k+2}-2]$, which are therefore not necessarily reduced. 
One could reduce the resulting element using the specialisation to GRPs
of~\cite[Algorithm 5]{CH3}, which shows how to do this for a general LWPFI. Chung and Hasan refer to
this process as {\em short coefficient reduction} (SCR), as opposed to full modular reduction. 
However, for ECC operations it is faster (and more secure) to simply ignore
this expansion and rely on a later modular multiplication to perform the reduction, as is
required when computing a point addition or doubling,
see~\cite{25519,bernstein2} and \S\ref{sec:secure}.

\subsubsection{Squaring}

When $t$ is single precision, the CVMA formulae do not have any common subexpressions,
as arises for ordinary integer residue squaring. In this case GRP squaring
is performed using Algorithm~\ref{modmul}.
If $t$ is multiprecision, then the components of a product $\vec{x} \cdot \vec{y}$ are computed as a sum 
of integer squares. In this case, one can eliminate common subexpressions 
to improve efficiency by nearly a factor of two (in the multiplication step).
On the other hand, when using Algorithm~\ref{interleaved} and its variants it may be difficult to eliminate common
subexpressions efficiently~\cite{koc}.

\subsubsection{Inversion and equality check}

Inversion seems difficult to perform efficiently in the GRP representation. If
$t$ were prime then it would be possible to use an analogue of the inversion/division
algorithm of~\cite{ozturk}, exploiting the cyclicity $t^{m+1} \equiv 1 \pmod{t^{m+1}-1}$. 
However, for our GRPs $t$ is even and greater than $2$. One can therefore opt to map back to $\Z/\Phi_{m+1}(t)\Z$, remaining in
the Montgomery domain, and perform inversion using the binary extended
Euclidean algorithm (see~\cite{knuth}, for example) and modular multplying by the precomputed value
$\psi(b^3 \bmod{\Phi_{m+1}(t)})$. Alternatively, for data-independent inversion, one can simply power by
$\Phi_{m+1}(t) - 2$, as do the authors of~~\cite{bernstein2}. 
Using projective coordinates can obviate the need for inversions altogether,
however for many protocols inversion is unavoidable and when it is avoidable,
in some scenarios such representations of points should be randomised after a point multiplication~\cite{projleak}.

Since our representation possesses redundancy, equality checking is naturally problematic. We therefore
opt to map back to $\Z/\Phi_{m+1}(t)\Z$ to check equality there --- as for
inversion --- while remaining in the Montgomery domain. For ECC equality checking is usually a one-time computation
per coordinate, and so again this operation does not greatly impinge upon efficiency.

\subsection{Side-channel secure ECC}
\label{sec:secure}

As we demonstrated in \S\ref{sec:stability}, by choosing $t$, $l$ and $m+1$ carefully, one can avoid the
need to compute any final additions or subtractions when performing a modular
reduction. This is an analogue to various results for ordinary Montgomery arithmetic~\cite{walter3,walter4,hachez}.
The lack of a conditional addition/subtraction averts threats such
as~\cite{walter2,sakai}, the latter of which applies directly to the NIST GMNs.
Our modular multiplication algorithm is thus control-flow invariant with no
data-dependent operations, making it immune to timing attacks~\cite{kocher1} and simple power
analysis (SPA).

In addition to making modular multiplications and squarings immune to timing
attacks and SPA, one can also ensure that the computation of an entire elliptic curve point addition or
doubling is also immune. To do so, one chooses a GRP with $t$ divisible
by a sufficiently high power of $2$, so that during the course of
an elliptic curve point operation, even if one ignores the coefficient expansion caused
by additions/subtractions, these do not overflow and the modular reductions ensure
the outputs are fully reduced elements. Note that this requires $b = 2^l \mid t$ to be
a few bits longer than the minimum $l$-values listed in
Tables~\ref{stableparameters1} and~\ref{stableparameters2}: for reasons of
space we do not include the analysis here.
By doing so, a point addition or doubling becomes an atomic
operation, where the sequence of arithmetic operations is entirely data-independent. In this case one only
needs point-multiplication-level defences against timing attacks and SPA, such
as the double-and-add-always algorithm due to Coron~\cite{coron},
or the use of Edwards curves, for which the addition formula can also be used for doubling~\cite{edwards}.
Hence, ECC over GRPs may be straight-line coded, which is beneficial for both efficiency
and security.

Lastly, our representation can also be made immune to differential power analysis (DPA)~\cite{kocher2}.
Observe that the embedding of $\Z/\Phi_{m+1}(t)\Z$ into $\Z/(t^{m+1}-1)\Z$ can be
randomised by adding to it a random multiple $r \cdot \Phi_{m+1}(t)$ for
$r \in \{0,\ldots,(t-1)-1\}$. While our embedding is an example of `operand
scaling'~\cite{walter,ozturk} which is used for faster reduction,
the addition of a multiple of the modulus within a redundant scaled representation 
also acts as a countermeasure to DPA --- such as Goubin's attack~~\cite{goubin} on the randomised projective
coordinates defence of Coron~~\cite{coron} --- as shown by Smart, Oswald and Page~\cite{danrandom}.
In particular, for multiprecision integer residues the authors show that this countermeasure thwarts DPA whenever
the scaling factor is longer than the longest string of ones or zeros in the
binary expansion of the initial modulus. For the NIST GMNs, this
countermeasure requires a large scaling factor, making the defence
inefficient and nullifying the benefits of using these moduli.
Applying the Smart-Oswald-Page rationale to GRPs, one sees that the scaling factor is $t-1$, while the longest string of ones or
zeros in the binary expansion of $\Phi_{m+1}(t)$ is $\lceil \log{t} \rceil -1$.
Since GRPs {\em already} use the larger ring, we acquire
this defence for almost negligible cost. In particular the
addition of a random multiple $r$ of $\Phi_{m+1}(t)$ to an element $\vec{x}$ has the form
$[x_m + r,x_{m-1}+r,\ldots,x_0 + r]$, which only requires $m+1$ additions.
Since DPA depends on the ability of an attacker to predict a specific bit in
the representation of a given field element (other than the upper excess zero bits
in each coefficient of GRP residues, which are the same for every field element), 
if the representation of points is randomised in this way prior to
every point multiplication, or even every modular multiplication, then DPA
should not be feasible.

%---------------------------------------------------------------------------------------------------
%---------------------------------------------------------------------------------------------------

\section{GRP Parameters}\label{sec:paramgen}

In this section we provide empirical data regarding the abundancy of
GRPs at various bitlengths relevant to ECC. We also specify parameters that are particularly suitable
for efficient implementation.

\subsection{Estimating the number of GRP parameters}\label{estimateparams}

As we saw from Tables~\ref{stableparameters1} and~\ref{stableparameters2}, for a
given prime $m+1$ and word size $w$, there is an upper bound on the length of
a GRP that may be represented. Table~\ref{paramcount} contains
estimates (or exact counts) for the number of GRPs which are in accordance with the GRP field and residue
representation set out in this work, for a word size $w = 64$ and where $q=2$
reductions suffice to ensure I/O modular multiplication stability.
The data was obtained as follows.

For a desired GRP $p$ of bitlength $\lceil \log{p} \rceil$, Table~\ref{stableparameters1} gives 
the minimum value of prime $m+1$ which is adequate to represent GRPs of this
size. The inequality~(\ref{constraint}) gives
$k_{max}$ which is the maximum bitlength of $t$ that is representable, 
while~(\ref{eq:lbound}) gives the minimum value $l$ required in order for $t = 2^l \cdot c$ 
to be I/O stable. We estimate $t_{max}$ simply as $2^{\lceil \log{p}\rceil /m}$, which implies
a maximum value for $c$ of $2^{\lceil \log{p}\rceil/m-l_{min}}$. Similarly for $p$ of
this precise bitlength, we estimate the minimum value of $c$ as $2^{(\lceil
  \log{p}\rceil-1)/m-l_{min}}$. We denote this interval by $I(c)$.
To estimate $P(\text{prime})$, which is the probability that a given
generalised repunit in our form is a GRP, we performed a linear search over $c
\in I(c)$, counting the first $1,000$ primes and simply dividing by the length of the
search. The estimated total number of GRPs satisfying our requirement that
$q=2$ is then given by $|I(c)| \cdot P(\text{prime})$.

\begin{table}\caption{Estimated GRP counts for $w=64$, $q=2$}\label{paramcount}
\begin{center}
\begin{tabular}{|c|c|c|c|c|c|c|c|c|}
\hline
$\lceil \log{p} \rceil$ & $m+1$ & $k_{max}$  & $\log{t_{max}}$ & $l_{min}$ & $|I(c)|$ & $P(\text{prime})$ & $\approx \#$GRPs\\
\hline
$600$ & $11$ & $60$ & $60.0$ & $34$ & $4.49 \times 10^6$ & $8.54 \times 10^{-3}$ & $38.4 \times 10^3$\\
$599$ & $11$ & $60$ & $59.9$ & $34$ & $4.19 \times 10^6$ & $9.05 \times 10^{-3}$ & $37.9 \times 10^3$\\
\hline
$512$ & $11$ & $60$ & $51.2$ & $30$ & $1.61 \times 10^5$ & $1.05 \times 10^{-2}$ & $1697$\\
$511$ & $11$ & $60$ & $51.1$ & $30$ & $1.51 \times 10^5$ & $1.06 \times 10^{-2}$ & $1591$\\
\hline
$384$ & $11$ & $60$ & $38.4$ & $24$ & $1448$ & $9.67 \times 10^{-3}$ & $14$\\
$383$ & $11$ & $60$ & $38.3$ & $24$ & $1352$ & $1.33 \times 10^{-2}$ & $18$\\
\hline
$360$ & $7$ & $60$ & $60.0$ & $34$ & $4.49 \times 10^6$ & $1.82 \times 10^{-2}$ & $81.7 \times 10^3$\\
$359$ & $7$ & $60$ & $59.9$ & $34$ & $4.19 \times 10^6$ & $1.77 \times 10^{-2}$ & $74.1 \times 10^3$\\
\hline
$256$ & $7$ & $60$ & $42.66$ & $25$ & $2.27 \times 10^4$ & $2.47 \times 10^{-2}$ & $561$\\
$255$ & $7$ & $60$ & $42.5$ & $25$ & $2.02 \times 10^4$ & $2.63 \times 10^{-2}$ & $531$\\
\hline
$244$ & $5$ & $61$& $61.0$ & $34$ & $2.14 \times 10^7$ & $1.68 \times 10^{-2}$ & $3.58 \times 10^5$\\
$243$ & $5$ & $61$& $60.75$ & $34$ & $1.80 \times 10^7$ & $1.72 \times 10^{-2}$ & $3.08 \times 10^5$\\
\hline
$224$ & $5$ & $56$ & $56.0$ & $31$ & $5.34 \times 10^6$ & $1.98 \times 10^{-2}$ & $1.06 \times 10^5$\\
$223$ & $5$ & $56$ & $55.75$ & $31$ & $4.49 \times 10^6$ & $1.88 \times 10^{-2}$ & $8.42 \times 10^4$\\
\hline
\end{tabular}
\end{center}
\end{table}

For each of $m+1 = 5,7$ and $11$, Table~\ref{paramcount} contains estimated counts
for the largest GRPs representable. It also contains estimates (or exact
counts) for the number of GRPs at the NIST GMN sizes $224,256$ and $384$. We also
consider bitlength $512$ rather than $521$, since this conjecturally gives
$256$-bit security, with the larger prime $2^{521}-1$ being nominated purely for
fast reduction. Observe that the number of available GRPs for a given $m+1$
decreases as the size of $p$, and hence $c$ decreases. The number available for
bitlengths $383$ and $384$ is particularly low. However, should this be a concern for a
particular application, one can
see from Table~\ref{stableparameters2} that by moving to GRPs for which $3$
reductions suffices, $|I(c)|$ becomes much larger ($3.71 \times 10^5$) and our
estimate of the number of GRPs becomes over $5,000$. 
On the other hand, since $384$ is not too far beyond the upper bound for the
size of GRP representable by $m+1 = 7$, it may be preferable to trade
$12$-bits of security for much improved performance, see \S\ref{sec:results}.
Similarly, 

\subsection{Hamming weight $2$ parameters}

As we showed in \S\ref{sec:redfast}, there are no suitable GRPs in the ECC
range for which $t=2^l$. Hence the next best type of GRP parameter $t$ will
have Hamming weight equal to $2$, where we allow $c$ to have the form
$2^{c'} + 1$ as well as $2^{c'}-1$ when there is sufficient slack in the representation,
since subtractions cost the same as additions. We list these GRPs in Table~\ref{fastparams}.
The final column indicates whether or not the given GRP allows for atomic side-channel secure point additions
and doublings, as per \S\ref{sec:secure}.
Note that for $m+1 = 5$ and $w=64$ we can not represent GRPs any larger than
$244$-bits, and are thus short of the conjectured $128$-bit ECC security level of $256$-bits. 
One can therefore either move up to $m+1=7$, which can represent GRPs of up to
$360$-bits, or one can opt to reduce security by a few bits, for better performance.
Indeed, in recent work K\"{a}sper argues that the NIST GMN prime $\text{P-224} = 2^{224} - 2^{96} + 1$
offers a satisfactory trade-off between security and efficiency, when used as
the basis of the elliptic curve Diffie-Hellman (ECDH) key exchange in the Transport
Layer Security (TLS) protocol~\cite{tls}. Bernstein has also implemented
arithmetic mod $\text{P-224}$~~\cite{p224}.
Yet another possibilty at this security level are the GFN primes 
$\Phi_8(2^{41}\cdot(2^{15}-1))$ and $\Phi_8(2^{50}\cdot(2^{6}-1))$, both of which have bitlength $224$, 
but experiments with such GFNs have not yet been carried out.

\begin{table}\caption{Approximately NIST-size fast GRPs for $w=64$, $q=2$}\label{fastparams}
\begin{center}
\begin{tabular}{|c|c|c|}
\hline
$\lceil \log{p} \rceil$ & GRP & S.C. Secure\\ 
\hline
$511$ & $\Phi_{11}(2^{42} \cdot (2^9 + 1))$ & Yes\\
\hline
$381$ & $\Phi_{11}(2^{34} \cdot (2^{4}-1))$ & Yes\\
$380$ & $\Phi_{11}(2^{34} \cdot (2^{4}+1))$ & Yes\\
\hline
$270$ & $\Phi_7(2^{34} \cdot (2^{11}-1))$ & Yes\\
$253$ & $\Phi_7(2^{27} \cdot (2^{15}+1))$ & Yes\\
$253$ & $\Phi_7(2^{37} \cdot (2^{5}+1))$ & Yes\\
$243$ & $\Phi_5(2^{59} \cdot (2^2-1))$ & No\\
\hline
$228$ & $\Phi_{5}(2^{54} \cdot (2^3-1))$ & Yes\\
$224$ & $\Phi_{5}(2^{31}\cdot(2^{25}-1))$ & No\\
$220$ & $\Phi_{5}(2^{52} \cdot (2^3-1))$ & Yes\\
\hline
\end{tabular}
\end{center}
\end{table}

\section{Implementation and Results}\label{sec:results}

In this section we provide details of our proof-of-concept implementation and
our results. We consider field multiplications only as this is the bottleneck for ECC point multiplication and hence 
gives an accurate indication of performance.

In terms of performance, the fastest implementations of ECC in the literature all feature cycle counts
for $256$-bit ECC point multiplication~\cite{25519,speed,hisil,longa,GLS,bernstein2}, except for K\"{a}sper's $\text{P-224}$ 
implementation~\cite{kasper}, with~\cite{25519,bernstein2,kasper} each being side-channel secure.
As it is difficult to get a fair comparison between our implementation and these, we opt to compare our modular 
multiplication performance with the \url{mp}$\F_q$ benchmarking system due to Gaudry and Thom\'{e}~\cite{speed}. 
This has been ported to OS-X 10.5.8 with minor changes and executed on a platform using an Intel Core 2 
Duo at 2.2Ghz. As stated in~\cite{bernstein2}, to date \url{mp}$\F_q$ gives only the fourth fastest implementation of 
ECDH, based on Bernstein's \url{curve25519}, which utilises a non-standard representation of residues mod $2^{255}-19$ and exploits
the floating-point unit of specific instruction-set architectures to great effect.
However, by comparing the basic multiplication cost on the target architecture, one can obtain a crude estimate of the relative
performance of our arithmetic with that of \url{curve25519}. 

Our implementation consists of two inline assembly operations targeted at the Core 2 processor. 
One accumulates the innermost sum of \url{line} \url{2} of Algorithm~\ref{modmul}, while the other performs a single 
instance of the reduction operation in \url{line} \url{4} of Algorithm~\ref{modmul}. 
Both use the 64-bit operations available on the Core 2 and the 
extended register set available in x86\_64. 
%Observe that since the size of the input's components are shorter for the
%second invocation of the reduction function than the first, it is necessary to
%implement two separate reduction functions. Also don't we need another two for general $c$? 
These assembly operations both use a mere $4$ of the $15$ available in the x86-64 instruction
set. This allows one to rely on normal C code to arrange these macros, and to handle data-storage. As  
a result the gcc compiler can generate all of the intermediate memory access instructions and schedule
the usage of the other $11$ registers available. This means that the same code can be reused
for any field supported by Algorithm~\ref{modmul} --- the only changes required are the parameter definitions.
To generate a particular instance of the family of algorithms we use a simple wrapper written in Python that 
arranges the sequence of these operations required for the particular parameter choice of $m+1$ and $t$. 

To emphasise the relative simplicity of our implementation, we use only 64-bit scalar operations on the processor, 
and allow the compiler to schedule most of the output instructions. As a result we reach a throughput of slightly 
less than one operation per cycle. In comparison the \url{mp}$\F_q$ implementation of \url{curve25519} uses SSE2 to reach a throughput of almost 
two operations per cycle (the theoretical maximum on the architecture). Although our implementation is less efficient 
(because we have spent less programmer time on the machine-dependent optimisation) the performance achieved is still 
higher. Scheduling a lower-level implementation on the processor would be an interesting challenge.

As explained in \S\ref{sec:reduction}, within the reduction algorithm we have a trade-off between the
number of GRPs available and performance. If one opts for a generic value of $c$ many GRPs are available, 
but the reduction involves a full \url{imulq} instruction with relatively high latency. If we specialise our choice  
of $c$ to very low Hamming weight then we can replace this instruction with a combination of  
shift and add instructions to improve performance. We have measured  
the performance of both implementations. To ensure a fair comparison we have merged our code 
into \url{mp}$\F_q$ so that all algorithms are being tested with 
the same timing code. This timer executes $10^6$ operations in the field, measuring the
elapsed time. The reported figures are the mean execution time for the operation.
Table~\ref{mont} contains cycle counts for Montgomery arithmetic at various bitlengths,
as well as the \url{curve25519} modular multiplication cycle count.
Table~\ref{mrcp} contains our results for GRP modular multiplication.

\begin{table}\caption{mp$\F_q$ cycle counts for curve25519 and Montgomery arithmetic}\label{mont}
\begin{center}
\begin{tabular}{|c|c|c|}
\hline
Algorithm & Size (bits) & Mult (cycles) \\
\hline
M. & $64$ & $30$ \\
M. & $128$ & $105$ \\
M. & $192$ & $195$ \\
\url{curve25519} & $255$ & $140$ \\
M. & $256$ & $280$ \\
M. & $320$ & $407$ \\
M. & $384$ & $563$ \\
M. & $448$ & $757$ \\
M. & $512$ & $981$ \\
\hline
\end{tabular}
\end{center}
\end{table}
\begin{table}\caption{Cycle counts for GRP arithmetic}\label{mrcp}
\begin{center}
\begin{tabular}{|c|c|c|}
\hline
Parameters & Max size (bits) & ModMul (cycles) \\
\hline
$m + 1\hspace{-2pt}=\hspace{-2pt}5$, $HW(c)\hspace{-2pt}=\hspace{-2pt}2$ & $244$ & $96$ \\
$m + 1\hspace{-2pt}=\hspace{-2pt}5$, general $c$ & $244$ & $112$ \\
$m + 1\hspace{-2pt}=\hspace{-2pt}7$, $HW(c)=2$ & $360$ & $165$ \\
$m + 1\hspace{-2pt}=\hspace{-2pt}7$, general $c$ & $360$ & $182$ \\
$m + 1\hspace{-2pt}=\hspace{-2pt}11$, general $c$ & $600$ & $340$ \\
\hline
\end{tabular}
\end{center}
\end{table}

As stated in \S\ref{estimateparams}, the closest size of field to \url{curve25519} that we can implement using $m+1 = 5$ 
is only $244$-bits. This small reduction in field size is  
compensated by an increase in performance, requiring only $80\%$ of the \url{curve25519} cycles per 
multiplication. Using the specialised reduction function for the $243$-bit GRP
$\Phi_{5}(2^{59}\cdot (2^2 -1))$, this figure improves to $69\%$.
Since the results for the first line of Table~\ref{mrcp} apply also to
Hamming weight $2$ GRPs smaller than $2^{243}$, we obtain the same
modular multiplication performance, while utilising the acquired slack in the
representation to ensure atomic point doublings/additions as per \S\ref{sec:secure}, in particular for the $228$-bit GRP 
$\Phi_{5}(2^{54}\cdot (2^3 -1))$. At bitlength $512$ with general $c$, compared to Montgomery multiplication, GRP multiplication
costs only $35\%$ as many cycles. At bitlength $600$, this proportion would naturally be even smaller, however 
at this size Karatsuba multiplication may be faster than schoolbook arithmetic.
We thus expect that point multiplications at $224$-bits and $512$-bits using GRPs to be competitive with the state-of-the-art
in the literature.

We freely admit that our proof-of concept implementation has not been
optimised, and therefore believe that one could obtain significantly better performance figures.
By comparing our arithmetic with the modular multiplication used in~\cite{bernstein2}, which is the benchmark for point multiplication
at the $128$-bit security level, one gains an idea of the potential performance of arithmetic mod $\Phi_{5}(2^{54}\cdot (2^3 -1))$ for example.
In~\cite{bernstein2}, residues are also represented by five $64$-bit words. Residue multiplication requires
$25$ \url{mul} instructions, as well as some \url{imul}, \url{add} and \url{adc} instructions. 
In comparison, to multiply $\vec{x}$ and $\vec{y}$ in our representation,
the CVMA formulae are as follows:
\begin{eqnarray}
\nonumber z_0 &=& (x_4 - x_1)(y_1 - y_4) + (x_3 - x_2)(y_2 - y_3),\\
\nonumber z_1 &=& (x_2 - x_4)(y_4 - y_2) + (x_1 - x_0)(y_0 - y_1),\\
\nonumber z_2 &=& (x_0 - x_2)(y_2 - y_0) + (x_4 - x_3)(y_3 - y_4),\\
\nonumber z_3 &=& (x_3 - x_0)(y_0 - y_3) + (x_2 - x_1)(y_1 - y_2),\\
\nonumber z_4 &=& (x_1 - x_3)(y_3 - y_1) + (x_0 - x_4)(y_4 - y_0),
\end{eqnarray}
requiring only $10$ \url{mul}, $25$ \url{add} and $5$ \url{adc} instructions. 
Since the respective reduction algorithms are quite similar with both requiring two rounds of
shifts, masks and additions, one expects the GRP modular multiplication to be considerably
faster, when optimsed. It is also possible that an optimised implementation of multiplication mod the $m + 1=7$ GRPs 
listed in Table~\ref{fastparams} may be faster than~\cite{bernstein2}, since it requires $21$ \url{mul} instructions, rather than $25$.
However, since this paper is predominantly expositional, we leave such optimisations as open research.

%----------------------------------------------------------------------------------------------------
%----------------------------------------------------------------------------------------------------

\section{Conclusion}
\label{sec:conclusion}

We have proposed efficient algorithms for performing arithmetic modulo a large family of primes, namely the generalised
repunit primes. The algorithms are simple to implement, are fast, are easily 
parallelisable, can be made side-channel secure, and all across a wide range of field sizes. 
The central contribution of this work is the development of the necessary theory, covering field and residue representation, 
as well as novel algorithms for performing efficient multiplication and reduction in these fields.
We have also presented proof-of-concept implementation results which provide an empirical comparison with 
other results in the literature, ensuring a fair comparison by reusing the same benchmarking procedure. 
Against Montgomery arithmetic we show an approximate three-fold increase in performance, and expect optimised implementations
of point multiplications using our proposed family to be competitive with the state-of-the-art in the literature.
We thus present a compelling argument in favour of a new approach to the secure and efficient implementation of ECC.

\section*{Acknowledgements}
The authors would like to thank Dan Page for making several very useful comments and suggestions, and the referees for their comments.

%------------------------------------------------------------------------------------------------
%------------------------------------------------------------------------------------------------

\bibliographystyle{amsplain}

\end{document}